\documentclass[11pt]{article}
\usepackage[hscale=0.8,vscale=0.8]{geometry}

\usepackage{amsmath}
\usepackage{amsfonts}
\usepackage{latexsym}
\usepackage{graphicx}
\usepackage{caption}
\usepackage{subcaption}
\usepackage{multicol}
\usepackage{multirow}
\usepackage{amssymb}
\usepackage{mathtools}
\usepackage{hhline}
\usepackage{amsthm}
\usepackage{float}
\usepackage{cite}
\usepackage{color}
\usepackage{chngcntr}
\usepackage{tikz}
\usepackage{xcolor}
\usepackage{array}
\usepackage{comment}
\usepackage{geometry}
\usepackage{bm}
\usepackage{epsfig}
\usepackage{algorithm}
\usepackage{algpseudocode}

\graphicspath{{./figs/}}

\usepackage[normalem]{ulem}
\usepackage[normalem]{ulem}

\newcommand{\ds}{\displaystyle}

\def\qon{{\quad\hbox{on}\quad}}
\def\qan{{\quad\hbox{and}\quad}}

\newcommand{\bk}{\mathbf{k}}
\newcommand{\bu}{\mathbf{u}}
\newcommand{\bv}{\mathbf{v}}

\newcommand{\0}{{\mathbf{0}}}

\newcommand{\bsi}{{\boldsymbol\sigma}}
\newcommand{\btau}{{\boldsymbol\tau}}

\def\bF{\mathbf{F}}

\def\bJ{\mathbf{J}}

\numberwithin{equation}{section}
\allowdisplaybreaks

\allowdisplaybreaks

\newtheorem{remark}{Remark}[section]

\title{Improving numerical accuracy for the viscous-plastic formulation of sea ice}

\author{
  {\sc Tongtong Li}\thanks{Department of Mathematics, Dartmouth College, Hanover, NH 03755, USA, email: {\tt \{tongtong.li@dartmouth.edu, annegelb@math.dartmouth.edu, yoonsang.lee@dartmouth.edu\} }. All authors are supported by ONR MURI \#N00014-20-1-2595. AG is also supported in part by NSF grants  DMS \#1502640 and DMS \#1912685, and  AFOSR grant \#FA9550-18-1-0316. YL is also supported by NSF grant DMS \#1912999.}	
\quad
    {\sc Anne Gelb}\footnotemark[1]~
\quad    
    {\sc Yoonsang Lee}\footnotemark[1]~}

\date{\today}

\begin{document}

\maketitle

\begin{abstract}

Accurate modeling of sea ice dynamics is  critical for predicting environmental variables and is important in applications such as navigating ice breaker ships.  Research for both modeling and simulating sea ice dynamics is ongoing, with the most widely accepted model based on the viscous-plastic (VP) formulation introduced by Hibler in 1979. Due to its highly nonlinear features, this model is intrinsically challenging for computational solvers. In particular, sea ice simulations often significantly differ from satellite observations. This study therefore focuses on improving the numerical accuracy of the VP sea ice model. Since the poor convergence observed in existing numerical simulations stems from the nonlinear nature of the VP formulation,  this investigation proposes using the celebrated weighted essentially non-oscillatory (WENO) scheme -- as opposed to the frequently  employed  centered difference (CD) scheme -- for the spatial derivatives in the VP sea ice model.   We then proceed to numerically demonstrate that WENO yields  higher-order convergence for smooth solutions, and that furthermore it is able to resolve the discontinuities in the sharp features of sea ice covers -- something that is not possible using CD methods. Finally, our proposed framework integrates a potential function method that utilizes the phase field method to naturally incorporates the physical restrictions of ice thickness and ice concentration in transport equations, resulting in  a modified transport equations which includes additional forcing terms. Our method does not require post-processing, thereby avoiding the possible introduction of discontinuities and corresponding negative impact on the solution behavior.  Numerical experiments are provided to demonstrate the efficacy of our new methodology.

\end{abstract}

\section{Introduction}
\label{sec:introduction}
Sea ice dynamics plays a vital role in understanding the ice cover in polar regions. Properly representing sea ice dynamics is crucial in predicting environmental variables and is important in a wide range of applications such as the navigation of ice breaker ships \cite{Parno19, Wang21}. The observations of the Arctic Ice Dynamics Joint Experiment (AIDJEX) significantly improved the modeling of the sea ice dynamics in the 1970s \cite{Hunke10}. Since then there has been an increased effort in modeling sea ice dynamics \cite{Hibler79, Hunke97, Hunke01, Girard11, Wilchinsky06, Dansereau16}. The VP sea ice model introduced by Hibler \cite{Hibler79} has become the most widely used approach for sea ice dynamics. The model consists of a nonlinear momentum equation and two transport equations, and is initially developed for meshes in the range of 100 kilometers. To solve the momentum equation at this spatial resolution, implicit time-stepping schemes are recommended \cite{Ip91} due to the nonlinear character of the momentum equation stemming from the viscous-plastic material law. Solvers such as Jacobian-free Newton-Krylov (JFNK) solver \cite{Lemieux10}, which is an improved approximation to the Jacobian matrix (typically used in other implicit numerical simulations), have been developed to improve the numerical efficiency for solving the VP model. To avoid implicit methods altogether, the Elastic-Viscous-Plastic (EVP) model is proposed by Hunke and Dukowicz \cite{Hunke97} to relax the stability condition for an explicit time-stepping scheme. 

With increasing mesh resolutions now available, it is becoming increasingly apparent that the numerical solutions to the sea ice model resulting from either the VP or EVP formulation are not well resolved, and indeed there is a significant discrepancy with obtainable observations \cite{Kwok08}. How much of this discrepancy is attributable to modeling error and how much to the numerical approximation error remains an open question \cite{Lemieux09}. In this study we focus on improving the numerical accuracy of the sea ice representation based on the one-dimensional VP model\footnote{We will discuss both the VP and EVP models. Since the EVP model is based on the VP model, henceforth, for ease of presentation, we will simply write VP when discussing the general model and only write EVP when needed for clarity.} which will be discussed in Section~\ref{sec:seaicemodel}.  

While there are many aspects of the VP model that  merit investigation, our focus here is on two different but related issues: (1) the numerical efficacy of the computational methods used for solving the model, which includes a study on both accuracy and convergence; and (2) ensuring that the computational method observes physical constraints on the ice thickness and concentration. Discussion on each issue as motivation for this investigation is provided below.

\subsection{Numerical efficacy of the VP sea ice model}
\label{subsec:VPmodelconvergence}

One goal of this investigation is to address the numerical simulation of the nonlinear VP sea ice model, specifically  concerning the accuracy, stability, and efficiency. While many efforts have been made to improve the computational {\em efficiency} of sea ice model solvers, e.g. \cite{Lemieux10, Lemieux12, Hunke97, Kimmritz15}, analysis of the corresponding {\em convergence} properties is lacking, and indeed many of these methods fail to converge\cite{Lemieux14}. We note that \cite{Lemieux14}, which proposes and implements an iterated IMplicit-EXplicit (IMEX) time integration to solve the coupled sea ice model monolithically, does provide an analysis of the temporal convergence of the numerical solution. There it is demonstrated that a combination of the second-order Runge-Kutta method for the explicit time integration and a second-order backward difference method for the implicit integration of the momentum equation yields an overall second-order accuracy in time of the numerical solution when compared to a reference solution obtained using a tiny time step ($1$ second). Spatial convergence is investigated in \cite{Williams18}, where it is shown that in the VP sea ice model, the simulated velocity field depends on the spatial resolution of the model and approaches the analytical solution as the spatial resolution is increased. However, the study is not quantified in terms of convergence rate. The method in \cite{Seinen17} adopts the Crank Nicholson time discretization and centered difference (CD) spatial discretization together with the JFNK solver. Second-order convergence is obtained using a synthetic model by simultaneously refining the spatial resolution and time step. It is not observed everywhere during the simulation test, however. 

To the best of our knowledge, convergence with respect to spatial resolution has not been well studied. In particular, no clear conclusion has been drawn in terms of spatial convergence. We adopt this as a starting point in our investigation to explore the related numerical properties. As we focus on spatial convergence, we choose the time step to be sufficiently small so that the time discretization error does not affect the convergence rate. This allows us to test  the VP formulation using an explicit time-stepping scheme and therefore avoid the error caused by the nonlinear solver needed for the implicit time-stepping scheme. Based on a constructed analytical solution with appropriately added forcing term to the governing equations, we test for convergence on the VP sea ice model using  both a second-order CD spatial discretization scheme as well as a third-order total variation diminishing (TVD) Runge-Kutta time integration scheme. 

Nowadays, with increasing mesh resolutions of up to $10$ km \cite{Lemieux14}, the general performance of the existing sea ice solvers is degrading, leading to a significant increase in numerical cost \cite{Hunke10}. Higher-order spatial discretization methods may be able to offset this problem. Due to the natural discontinuity feature of ice concentration and ice thickness, traditional higher-order finite difference schemes typically have spurious oscillations near discontinuities (the Gibbs phenomenon), which may pollute smooth regions and even lead to instability, causing blowups of the schemes \cite{Shu20}. The weighted essentially non-oscillatory (WENO) method \cite{Liu94} is designed to achieve higher-order accuracy in smooth regions while sharply resolving discontinuities in an essentially non-oscillatory fashion.  This study verifies that these desirable properties hold when implementing WENO for the sea ice model. 

\subsection{Ice thickness and ice concentration}
\label{subsec:icethickness_Intro}

In most sea ice models,  including the VP model, two idealized thickness levels, namely thick and thin, are often adopted to approximately characterize ice thickness in a relatively simple form. The two variables used to keep track of these levels are ice thickness, which is equivalent to the mass of ice in any grid cell, and ice concentration, which is defined as the fraction of the grid cell area covered by thick ice. 

One issue is how the constraints on these two variables are imposed. In particular, assuming continuity of the ice thickness and the ice velocity, the ice thickness should remain non-negative (see \cite[Theorem 3.10]{Kuzmin10}). The non-negativity of ice concentration is similarly guaranteed. As will be demonstrated, preserving the non-negativity in both parameters is an important consideration for choosing a numerical solver. Moreover, although not explicitly providing a method to guarantee the upper bound of ice concentration to be $1$, such a constraint is described in Hibler \cite{Hibler79}  to be equivalent to adding a mechanical sink term in the model, which is turned on when the ice concentration reaches $1$ to prevent its further growth.

As far as we know, the numerical methods used to solve the model or impose the constraints have not been rigorously analyzed or compared. As already mentioned, the original work in \cite{Hibler79} does not explicitly discuss how these constraints may be incorporated into the numerical implementation of the model. Since then some investigations have explicitly provided approaches both for numerical simulation of the model as well as for imposing these constraints. For example, Mehlmann \cite{Mehlmann19} uses a finite element framework to solve the sea ice model and imposes restrictions on the trial spaces of ice thickness and ice concentration through a projection of the solution. Lipscomb and Hunke \cite{Lipscomb04} adopt an incremental remapping scheme for sea ice transport. This is a Lagrangian approach that preserves the monotonicity by Van Leer limiting. That is, the gradients are reduced when necessary to ensure that the values in the reconstructed fields stay inside the range of the mean values in the cell and its neighbors. All of these mentioned approaches impose the model constraints through a post-processing procedure, which may, unfortunately, introduce discontinuity into the numerical solution, which affects the accuracy and might ultimately impact stability so that the solution does not converge. 

Therefore, the other goal of this investigation is to develop a numerical approach that more intuitively imposes the model constraints without any post-processing procedure. To accomplish this task, we propose using the potential function method motivated by the analogous approach of using a double-well potential function in what is commonly referred to as the phase field method \cite{Fix83, Langer86}, which is designed to solve interface problems by treating the interface as an object with finite thickness. Our proposed method, described in Section \ref{sec:method}, offers a simple but elegant way to incorporate additional restrictions into the model and could be generalized in various settings. It further yields a modified transport model with the extra forcing terms coming from the potential energy function, which is consistent with Hibler's statement on how to include the mechanical sink term in the model. 

The rest of the paper is organized as follows. In Section~\ref{sec:seaicemodel} we describe the sea ice model. Section~\ref{sec:solver} provides a brief overview of standard numerical solvers for the sea ice model, focusing on the JFNK solver and EVP solver. We describe both the WENO scheme and the potential function method in Section~\ref{sec:method}. In Section~\ref{sec:numerical} we conduct some numerical experiments and compare the performances of the WENO method with the more typically employed CD scheme. We also provide some numerical illustrations depicting the use of the potential function method. We make some concluding remarks in Section~\ref{sec:conclusion}.

\section{Sea ice dynamics model}
\label{sec:seaicemodel}
We begin by describing the two-dimensional VP model introduced by Hibler \cite{Hibler79} for the simulation of sea ice circulation and thickness. Although sea ice dynamics occurs in a three-dimensional space, the vertical scale of ${\mathcal O}(m)$ is much smaller than the horizontal scale of ${\mathcal O}(1000\ m)$, so the motion of sea ice is usually described in two dimensions. The VP sea ice model comprises of a momentum equation and two transport equations that describe the balance laws and is given by
\begin{subequations}
\label{eq:2Dmodel}
\begin{equation}
\ds m \frac{D \bu}{Dt}= m(\frac{\partial \bu}{\partial t} +  \bu \cdot \nabla \bu) = \nabla \cdot \bsi-m f \, \bk \times \bu + \btau_a-\btau_w - m g \nabla H_d, \label{eq: momentum u 2D}
\end{equation}
\begin{equation}
\ds \frac{\partial h}{\partial t}+ \nabla \cdot (\bu \, h)=S_h, \label{eq: transport h 2D}
\end{equation}
\begin{equation}
\ds \frac{\partial A}{\partial t}+ \nabla \cdot (\bu \, A)=S_A. \label{eq: transport A 2D}
\end{equation}
\end{subequations}
Here $\bu$ is the two-dimensional ice velocity, $h$ is the mean ice thickness, and $A$ is the ice concentration. The ice mass per unit area $m$ is given by $\rho h$, where $\rho$ is the sea ice density. The internal forces are modeled by $\nabla \cdot \bsi$ where $\bsi$ is the internal ice stress. The external forces comprise of the Coriolis force, forces due to air and water stress $\btau_a$ and $\btau_w$, and the force to the surface height. The other parameters include $f$, the Coriolis parameter, $\bk$, a unit vector perpendicular to the horizontal plane, $g$ the acceleration due to gravity, and $H_d$ the sea surface dynamic height. Finally, $S_h$ and $S_A$ are the thermodynamic source or sink terms.
We note that the advection term $\bu \cdot \nabla \bu$ of ice momentum can be neglected due to scaling properties \cite{Zhang91}. Furthermore, the thermodynamic terms are set to zero in the simulations as we concentrate on dynamic effects.

To better understand and analyze how well different computational approaches are suited to sea ice dynamics, we focus on a simplified one-dimensional sea ice model in this study, which is given by
\begin{subequations}
\label{eq:1DVPmodel}
\begin{equation}
{\ds \rho h \frac{\partial u}{\partial t}-\tau_a+\tau_w-\frac{\partial \sigma}{\partial x}=0,} \label{eq: momentum u 1D}
\end{equation}
\begin{equation}
\ds \frac{\partial h}{\partial t}+ \frac{\partial}{\partial x}(u \, h)=0, \label{eq: transport h 1D}
\end{equation}
\begin{equation}
\ds \frac{\partial A}{\partial t}+ \frac{\partial}{\partial x}(u \, A)=0. \label{eq: transport A 1D}
\end{equation}
\end{subequations}
$u$ is the one-dimensional sea ice velocity, and $\sigma$ is the internal stress corresponding to $\bsi_{xx}$ of the 2D model. In this model, the Coriolis force and sea surface tile are set to zero as the external forces act only in one direction on a static ocean slab \cite{Lipscomb07, Williams18}. The air and water stress terms, $\tau_a$ and $\tau_w$ respectively, are determined from  the nonlinear boundary layer theories and the quadratic drag formulas used in the model \cite{McPhee75}:
\begin{gather}
\ds \tau_a = \rho_a C_{da} \vert u_a \vert u_a, \label{eq: tau a}\\[1ex]
\ds \tau_w = \rho_w C_{dw} \sqrt{u^2+\epsilon_1} u, \label{eq: tau w}
\end{gather}
where $\rho_a$ and $\rho_w$ are the air and water densities, $C_{da}$ and $C_{dw}$ are the air and water drag coefficients, $u_a$ is the surface wind, and $\epsilon_1$ is a very small value ($10^{-10} \text{ m}^2/\text{s}^2$) introduced for numerical stability. Here the sea ice drift speed is neglected in the air drag formulation as it is much slower than the wind speed. The water under the ice is assumed to be at rest, leading to the absence of the water velocity in the water drag formulation.

We now describe the rheology term modeling the ice interaction, a viscous-plastic constitutive law relating the stresses and the strain rates. Due to the dimension reduction, all other components $\bsi_{xy}$, $\bsi_{yx}$ and $\bsi_{yy}$ vanish, and therefore the divergence of the stress tensor in \eqref{eq: momentum u 1D} is reduced to
\begin{equation}
\ds \frac{\partial \sigma}{\partial x}=\frac{\partial}{\partial x} \left[ (\eta+\zeta)\frac{\partial u}{\partial x}-\frac{P}{2} \right], \label{eq: div sigma}
\end{equation}
where
\begin{equation}
\ds \eta=\zeta e^{-2} \qan \zeta = \frac{P}{2 \Delta} \label{eq: viscosities}
\end{equation}
are the bulk and shear viscosities modeled by a normal flow rule in the plastic state and are chosen as constant values in the viscous regime. Here $e$ is the eccentricity of the elliptical yield curve, and $\Delta$ in one dimension is obtained as
\begin{equation}
\ds \Delta=\left[(1+e^{-2})\left( (\frac{\partial u}{\partial x})^2 + \epsilon_2 \right)\right]^{1/2}, \label{eq: Delta}
\end{equation}
with $\epsilon_2=10^{-22} \text{ s}^{-2}$ as another small parameter introduced for numerical stability purposes.

The original viscous-plastic formulation \cite{Hibler79} realizes the viscous coefficients by capping them at some maximum values, leading to a rheology term that is not continuously differentiable with respect to velocity. To obtain a smooth formulation of the viscous coefficients, we follow Lemieux and Tremblay \cite{Lemieux09} and replace the expression of $\zeta$ by the hyperbolic tangent function
\begin{equation}
\ds \zeta = \frac{P}{2 \Delta_{\min}} \tanh (\frac{\Delta_{\min}}{\Delta}), \label{eq: zeta}
\end{equation}
with $\Delta_{\min}=2 \times 10^{-9} \text{ s}^{-1}$ in accordance with the $\zeta_{\max}$ definition in \cite{Hibler79}.

The ice strength $P$ is expressed as
\begin{equation}
\ds P=P^{\star} h \exp [-C(1-A)], \label{eq: pressure}
\end{equation}
where $P^{\star}$ and $C$ are the strength and concentration parameters.

Finally, Table \ref{table: parameter} provides the values for all of the physical parameters which we will use in our numerical experiments. These parameter values are typically used in the VP sea ice model \cite{Lemieux12, Lemieux14, Williams17}.

\begingroup
\def\arraystretch{1.1}
\begin{table}[ht!]
\begin{center}
\begin{tabular}{l l l }
\hline
Symbol    & Definition    & Value  \\ \hline
$\rho$ & Sea ice density & $900 \text{ kg}/\text{m}^3$\\ 
$\rho_a$ & Air density & $1.3 \text{ kg}/\text{m}^3$\\ 
$\rho_w$ & Water density & $1026 \text{ kg}/\text{m}^3$\\ 
$C_{da}$ & Air drag coefficient & $1.2 \times 10^{-3}$\\
$C_{dw}$ & Water drag coefficient & $5.5 \times 10^{-3}$ \\
$P^{\star}$ & Ice strength parameter & $27.5 \times 10^3 \text{ N}/\text{m}^2$\\
$C$ & Ice concentration parameter & 20 \\
$e$ & Ellipse ratio & $2$ \\ \hline
\end{tabular}
\end{center}
\caption{Physical parameters used in the VP sea ice model.}
\label{table: parameter}
\end{table}
\endgroup

\section{Numerical solvers}
\label{sec:solver}
Time splitting methods are standard for solving the coupled sea ice system \eqref{eq: momentum u 1D}--\eqref{eq: transport A 1D} \cite{Lemieux14}, and in general, they are widely used to cope with the complex coupled system, e.g., in \cite{Hibler79, Hunke01, Lemieux14, Kimmritz15}. The basic idea is to decouple the momentum equation \eqref{eq: momentum u 1D} from the transport equations and solve it first, and then use the updated momentum to solve the transport equations, \eqref{eq: transport h 1D} and \eqref{eq: transport A 1D}, together. The main difficulty here lies in the momentum equation due to the highly nonlinear feature of the viscous-plastic rheology.

To apply an explicit time-stepping scheme to the momentum equation, numerical stability dictates a time step on the order of $1$ second for a $100$ km grid resolution \cite{Ip91}, or equivalently $1/100$ second for a $10$ km resolution grid, which is a typical spatial resolution for earth system models. Because of this very restrictive time step, it is recommended in \cite{Ip91} to use implicit time-stepping for the momentum equation. Implicit time-stepping requires the use of iterative methods which are notoriously difficult for nonlinear problems, however. To alleviate this issue, a Picard solver designed to repeatedly solve simple linear systems was proposed in  \cite{Zhang91}. Further investigation in \cite{Lemieux09} demonstrated the impractical slow convergence of the Picard solver which ultimately motivated the development of an inexact Newton method, realized as the JFNK solver, in \cite{Lemieux10}. 

On the flip side, in order to entirely avoid implicit methods, the EVP model proposed by Hunke and Dukowicz \cite{Hunke97} and then further modified by Hunke in \cite{Hunke01} adds  an artificial elastic term  to the viscous-plastic constitutive equation, thereby relaxing the stability condition for an explicit time-stepping scheme. The basic idea of the EVP model is to approximate the VP solution by damping the resulting artificial elastic waves via subcycling \cite{Lemieux12}. 

Below we provide a brief overview of the JFNK and EVP solvers focusing on the one-dimensional case.

\subsection{The Jacobian-free Newton-Krylov (JFNK) solver}
\label{sec:JFNK}
We illustrate the procedure of the JFNK implementation with a backward Euler time integration scheme for the momentum equation \eqref{eq: momentum u 1D}. The time-discretized one-dimensional momentum equation is written as
\begin{equation}
\ds \rho h^{n-1} \frac{ u^{n}-u^{n-1}}{\Delta t}=\tau_a^{n}-\tau_w(u^n)+\frac{\partial \sigma(u^n,h^{n-1},A^{n-1}) }{\partial x},  \label{eq: momentum u BE}
\end{equation}
where the superscript $n$ denotes the current time level.  The numerical solution at the previous time level $n-1$ for \eqref{eq:1DVPmodel} is known. Let ${\bf u}^n = \{u_j\}_{j = 1}^N$ denote the approximation of $u^n$ obtained by some finite difference spatial discretization technique at each grid point $x_j$, $j = 1,\dots N$.\footnote{To avoid cumbersome notation, for the rest of this paper we will use ${\bf u}$ to depict the vectorized solution of \eqref{eq:1DVPmodel} (and not the continuous solution \eqref{eq:2Dmodel}).}  The spatial discretization scheme will be specified on both non-staggered and staggered grids in Section~\ref{sec:method}. For now we generically define the solution on $N$ grid points. At current time level $n$, we therefore seek a  solution to 
$$ \bF({\bf u}^n)=\0,$$
where $\bF({\bf u}^n)$ is the difference between the right- and left-hand sides of \eqref{eq: momentum u BE} following spatial discretization. 

Since we are focusing on a single time step, we can simplify the notation by dropping the superscript $n$, so that we seek the solution ${\bf u}={\bf u}^{n}$. Using the velocity solution at the previous time level as the initial value ${\bf u}^{(0)}$, we  iteratively solve a sequence of linearized systems to consecutively obtain ${\bf u}^{(1)}, {\bf u}^{(2)}, \cdots, {\bf u}^{(k)}, \cdots$ until some stopping criterion is satisfied. Algorithm \ref{alg: JFNK} summarizes the iterative technique. More details can be found in \cite{Lemieux08, Lemieux10, Auclair17}.
\begin{algorithm}[h!]
\caption{JFNK solver}\label{alg: JFNK}
\begin{algorithmic}
\State Start with an initial iterate ${\bf u}^{(0)}$ and calculate $\Vert \bF({\bf u}^{(0)})\Vert$, here $\Vert \cdot \Vert$ is the L2-norm.
\For{$k=1$ to $k_{\max}=150$}
\State Solve $\bF({\bf u}^{(k-1)})+\bJ({\bf u}^{(k-1)}) \delta {\bf u}^{(k)}=\0$ for $\delta {\bf u}^{(k)}$, where $\ds \bJ({\bf u}^{k-1}) \bv \sim \frac{\bF({\bf u}^{(k-1)}+\epsilon \bv)-\bF({\bf u}^{(k-1)})}{\epsilon}$ and $\epsilon=10^{-7}$.
\State Set ${\bf u}^{(k)}={\bf u}^{(k-1)}+\lambda \, \delta {\bf u}^{(k)}$, where $\ds \lambda=\left[1, \, \frac{1}{2}, \, \frac{1}{4}, \, \frac{1}{8} \right]$ is successively reduced until $\Vert \bF({\bf u}^{(k)})\Vert <\Vert \bF({\bf u}^{(k-1)})\Vert$ or until $\ds \lambda=\frac{1}{8}$.
\State Stop if $\Vert \bF({\bf u}^{(k)})\Vert < \gamma_{nl}\Vert \bF({\bf u}^{0})\Vert$ with $\gamma_{nl}=10^{-6}$.
\EndFor
\end{algorithmic}
\end{algorithm}


Observe that Algorithm \ref{alg: JFNK} is an inexact Newton's method as it approximates the Jacobian $\bJ$ by a first-order Taylor series expansion. The linear system of equations is, in general, solved by the preconditioned FGMRES method \cite{Saad93}, which is a Krylov subspace method. This method is introduced as a matrix-free approach because forming and storing the Jacobian matrix is prohibitively expensive in CPU time and storage \cite{Lemieux10}. In the one-dimensional case, we are able to obtain the matrix representation of $\bJ$ by applying it to the basis vectors. We can then solve the linear system using a  direct solver as the computational cost and efficiency in the one-dimensional case are not causes for concern.

\subsection{The Elastic-Viscous-Plastic (EVP) solver}
\label{sec:EVP}

In the EVP solver, the velocity at time level $n$ is obtained by {\em explicitly} solving the momentum equation from the previous time level $n-1$. In particular, the constitutive law was rewritten by Hunke and Dukowicz \cite{Hunke97} to include a time dependence on the stress tensor. The velocity is then solved together with stress during subcycling.
In the one-dimensional case \cite{Williams17}, the stress-strain relationship
\begin{equation}
\ds \sigma= (\eta+\zeta)\frac{\partial u}{\partial x}-\frac{P}{2}  \label{eq: sigma}
\end{equation}
is equivalent to
\begin{equation}
\ds \frac{\sigma}{ \eta+\zeta}+\frac{P}{2(\eta+\zeta)}=\frac{\partial u}{\partial x}.  \label{eq: du}
\end{equation}
By adding an artificial elastic strain with an elastic parameter $E$, we obtain
\begin{equation}
\ds \frac{1}{E}\frac{\partial \sigma}{\partial t} + \frac{\sigma}{ \eta+\zeta}+\frac{P}{2(\eta+\zeta)}= \frac{\partial u}{\partial x}.  \label{eq: dt sig}
\end{equation}
In the original version of the EVP model \cite{Hunke97}, the viscosities $\eta$ and $\zeta$ were held fixed throughout the subcycling procedure. However, because the viscosities were not regularly updated, such linearization of the internal stress term caused the computed principal stress states to lie outside the elliptical yield curve \cite{Hunke01}. To address this issue, Hunke \cite{Hunke01} proposed to include the viscosities within the subcycling, while simultaneously changing the definition of the elastic parameter $E$ to maintain the computational efficiency. Specifically, with $E$ defined in terms of a damping timescale for elastic waves and $T$ according to the equation $\ds E=\frac{\zeta}{T}$, \eqref{eq: dt sig} can be rewritten as
\begin{equation}
\ds \frac{\partial \sigma}{\partial t} + \frac{\sigma}{ (1+e^{-2})T}+\frac{P}{2(1+e^{-2})T}=\frac{\zeta}{T} \frac{\partial u}{\partial x}.  \label{eq: dt sig 2}
\end{equation}
The subcycling solution is advanced iteratively with subcycling time step $\Delta t_e$. This approach yields the time evolution of stress as a function of the velocity from the previous iterate  according to
\begin{equation}
\frac{\sigma^{s}-\sigma^{s-1}}{\Delta t_e} +\frac{\sigma^{s}}{ (1+e^{-2})T}+\frac{P^{n-1}}{2(1+e^{-2})T}=\frac{\zeta^{s-1}}{T}\frac{\partial u^{s-1}}{\partial x}, \label{eq: EVP sig}
\end{equation}
where the subcycling iterate is denoted  with the superscript $s$. 
With the newly calculated stress in \eqref{eq: EVP sig}, the velocity is subcycled according to
\begin{equation}
\rho h^{n-1} \frac{u^{s}-u^{s-1}}{\Delta t_e} = \tau_a^{s}-\tau_w^{s-1}+\frac{\partial \sigma^{s}}{\partial x}.   \label{eq: EVP u}
\end{equation}
Observe that in the EVP solver, the damping timescale $T$ is a tuning parameter satisfying $\Delta t_e < T < \Delta t$, which is in general set to be $T=0.36\Delta t$ following the documentation of the CICE model \cite{Hunke15}. In addition, we denote the number of subcycles by $N_{sub}$, satisfying $N_{sub} \times \Delta t_e= \Delta t$.

\begin{remark}
We note that neither the JFNK nor the EVP solver has entirely resolved the convergence issue. In particular, the JFNK solver is not robust, as it was demonstrated in \cite{Lemieux10} that the failure rate for the JFNK solver increases as the grid is refined. On the other hand, it was shown in \cite{Lemieux12} that the EVP approximate solution has notable differences from the reference solution, which becomes relatively more distinct with finer resolution.  Furthermore, neither solver provides an explicit way to treat the out-of-range issues for either ice thickness or ice concentration. Hence it is possible to obtain unrealistic physical values for either or both when applying the solvers directly to the sea ice model. Thus we are motivated to address these issues in the 1D case so that we are better able to subsequently solve the more complicated two-dimensional version in \eqref{eq:2Dmodel}.
\end{remark}

\section{Proposed numerical methods}
\label{sec:method}
Motivated by the above discussion, in this section we propose an approach to help mitigate the limitations of existing solvers for both the VP and EVP sea ice models.  We first discuss the WENO method \cite{Liu94} to advocate the use of higher-order methods for improving numerical accuracy and efficiency.  We then describe the potential function method as a means to incorporate the physical restrictions of the ice thickness and concentration on top of the existing numerical methods. This will help to alleviate the out-of-range issues. 

\subsection{Weighted essentially non-oscillatory (WENO) scheme}
\label{sec:WENO}

A main goal of this investigation is to demonstrate the advantages of using higher-order methods to solve the sea ice model. We use the WENO method \cite{Liu94} as a prototype for two reasons.  First, WENO is designed to have a higher-order convergence rate for smooth solutions than standard three point stencil CD schemes (which are second order), and second, WENO is able to maintain stable, non-oscillatory, and sharp discontinuity transitions so that it is suitable for sea ice with natural discontinuity feature of thickness and concentration. We verify that both of these properties hold in our numerical simulation of sea ice cover with and without  sharp features.  

We use the method of lines for time integration in each numerical test. To ensure stability and maintain the accuracy obtained in the spatial derivative approximation, we use the third-order TVD Runge-Kutta (TVRK3) \cite{Shu88} when employing both WENO and CD schemes. We then compare the numerical convergence properties of the results obtained using both methods. We implement the same time integration scheme to ensure that we are only comparing the spatial discretization performances for each scheme and not evaluating the time integration methods. Indeed, using an implicit time integration (backward Euler) for the momentum equation and explicit time integration (forward Euler) for the transport equations is common \cite{Lemieux12, Lemieux14, Auclair17}. Our numerical example in Section \ref{sec:numtestpotential} provides a case study for the mixed time integration approach. It is also, of course, possible to use a different spatial discretization for each  equation in the system. As we did not observe any advantage in this approach, to reduce complexity, we did not do this.


In the WENO case we construct a non-staggered grid  so that all variables are defined at the center of each grid cell. That is, we seek the solution at the $j  =1,\dots, N$ midpoints of each cell, $\ds x_{j-\frac{1}{2}} = x_j -\frac{1}{2}\Delta x$, with $x_j = j\Delta x$ and $\ds\Delta x = \frac{L}{N}$ where $L$ is the domain length. Following \cite{Liu11} for nonlinear degenerate parabolic equations, we use higher-order finite differencing for \eqref{eq: momentum u 1D}. In particular, $\ds \frac{\partial u}{\partial x}$ and $\ds \frac{\partial \sigma}{\partial x}$ are discretized using the fifth order finite difference WENO method for conservation laws \cite{Jiang96} based on the left-biased stencil and right-biased stencil, respectively.  The fifth order WENO method for conservation laws \cite{Jiang96} is also used to solve the transport equations \eqref{eq: transport h 1D} and \eqref{eq: transport A 1D}. 

On the other hand, in the CD case, we construct a one-dimensional version of the staggered Arakawa C-grid \cite{Arakawa77}, where the velocity $u$ is defined on vertices, $u_0, \cdots, u_N$, and the traces $h$ and $A$ are defined at the center of each grid cell, $h_{\frac{1}{2}}, \cdots, h_{N-\frac{1}{2}}$ and $A_{\frac{1}{2}}, \cdots, A_{N-\frac{1}{2}}$ respectively. Correspondingly, the stress $\sigma$, the viscosities $\eta$ and $\zeta$ and the ice strength $P$ are also defined at the center of each grid cell. To solve for the velocity $u$ in the momentum equation \eqref{eq: momentum u 1D}, we take $\ds h_i=\frac{1}{2}(h_{i+\frac{1}{2}}+h_{i-\frac{1}{2}})$ for $i=1, \cdots, N-1$. We then approximate $\ds \frac{\partial u}{\partial x}$ at each cell center as
$$ \{du\}_{i+\frac{1}{2}} = \frac{u_{i+1}-u_i}{\Delta x},$$
so that $\ds \sigma=(\eta+\zeta)\frac{\partial u}{\partial x}-\frac{P}{2}$ is defined at $x_{i+\frac{1}{2}}$.  This leads to the approximation $\ds \frac{\partial \sigma}{\partial x}$ at each vertex given by
$$ \{d\sigma\}_{i} = \frac{\sigma_{i+\frac{1}{2}}-\sigma_{i-\frac{1}{2}}}{\Delta x}.$$
Similarly, for the transport equation \eqref{eq: transport h 1D} of $h$, we approximate $\ds \frac{\partial}{\partial x}(u h)$ at each cell center as
$$\{d(uh)\}_{i+\frac{1}{2}}=\frac{(uh)_{i+1}-(uh)_i}{\Delta x}.$$
We similarly solve  $A$  using the same spatial discretization for the  transport equation \eqref{eq: transport A 1D}. 
Finally, we note that although it is not a realistic assumption, we assume periodic boundary conditions to avoid errors introduced by boundary approximation. That said, the complexities introduced at the boundaries are not fundamentally different from other PDE models.

In terms of time integration, we consider the ordinary differential equation
$$ \frac{du}{dt}= L(u),$$
where $L(u)$ is a discretization of the spatial operator. The TVRK3 method advances the current solution $u^n$ to the next time level $u^{n+1}$ according to Algorithm \ref{alg: TVRK3}.

\begin{algorithm}[h!]
\caption{TVRK3 time integration}\label{alg: TVRK3}
\begin{algorithmic}
\State Start with initial value  ${u}^{n}$ and calculate $L(u^n)$.  To calculate the solution at the next time step $n+1$:
\State $ u^{(1)}=u^{n}+ \Delta t \, L(u^{n})$,
\medskip
\State $ u^{(2)}=\frac{3}{4}u^{n} + \frac{1}{4}u^{(1)} + \frac{1}{4}\Delta t \, L(u^{(1)})$,
\medskip
\State $ u^{n+1}=\frac{1}{3}u^{n} + \frac{2}{3}u^{(2)}+\frac{2}{3}\Delta t \, L(u^{(2)}).$
\end{algorithmic}
\end{algorithm}

\subsection{Potential function method}
\label{sec:potential}

Due to its physical interpretation, the variable ice thickness $h$ in the sea ice model must remain non-negative. Similarly, the ice concentration value $A$ must be between 0 and 1. It is crucial for the numerical methods to preserve both of these properties to ensure that a meaningful solution is obtained. Motivated by the double-well potential function in the phase field method \cite{Fix83, Langer86}, where the resulting equation is limited to a particular set of prescribed values due to the local minima of the potential function, we develop the potential function method here to impose the corresponding restrictions of ice thickness and ice concentration.

We begin by illustrating the potential function method on the transport equation of ice concentration $A$ and note that the case for ice concentration $h$ similarly follows. First, to  restrict ice concentration $A$ so that $0 \leq A \leq 1$, we define a potential function in a piecewise manner as
\begin{gather}\label{eq: potential}
\ds f(A)=
\begin{cases}
\gamma_1 f_1(A), & \text{if} \quad A<0,\\
0, & \text{if}\quad 0\leq A \leq 1, \\
\gamma_2 f_2(A), & \text{if}\quad A > 1,
\end{cases} 
\end{gather}
where $f_1 > 0$ and $f_2 > 0$ for all $A$, and $\gamma_1 > 0$ and $\gamma_2 > 0$ are parameters chosen so that $f$ has minima on $[0,1]$. For example, if both $f_1$ and $f_2$ are linear functions, a particular form of $f$ might be 
\begin{gather}\label{eq: potential linear}
\ds f(A)=
\begin{cases}
-\gamma_1 A, & \text{if} \quad A<0,\\
0, & \text{if}\quad 0\leq A \leq 1, \\
\gamma_2 (A-1), & \text{if}\quad A > 1.
\end{cases} 
\end{gather}
The transport equation \eqref{eq: transport A 1D} is then modified by adding a forcing term given by the gradient of the potential. This has the effect of the ice concentration experiencing a gradient force that tracks down to the physical range $[0,1]$.  The resulting equation is given by
\begin{equation}\label{eq: transport A potential}
\ds \frac{\partial A}{\partial t}+ \frac{\partial}{\partial x} (u \, A)=-f'(A).    
\end{equation}
Observe that for the piecewise linear case, the forcing term  $-f'(A)$ is piecewise constant, meaning that the force transition  is not continuous.  To enable a more desirable smooth transition for this term we will instead choose both $f_1$ and $f_2$ to be quadratic and define $f$  as 
\begin{gather}\label{eq: potential quadratic}
\ds f(A)=
\begin{cases}
\gamma_1 A^2, & \text{if} \quad A<0,\\
0, & \text{if}\quad 0\leq A \leq 1, \\
\gamma_2 (A-1)^2, & \text{if}\quad A > 1.
\end{cases}
\end{gather}
The corresponding forcing term is now given by
\begin{gather}
\label{eq:forcing}
\ds f'(A)=
\begin{cases}
2\gamma_1 A, & \text{if} \quad A<0,\\
0, & \text{if}\quad 0\leq A \leq 1, \\
2\gamma_2 (A-1), & \text{if}\quad A > 1,
\end{cases}
\end{gather}
which is clearly continuous.
%


\subsubsection{Determining parameters $\gamma_1$ and $\gamma_2$}
We now must determine parameters $\gamma_1$ and $\gamma_2$ in \eqref{eq: potential linear} that ensure $A$ will stay in range, that is $0 \le A \le 1$.  To this end, we first prescribe a Lagrangian representation of the ice concentration field $A$, which we will denote as
$$B(t) = A(x(t),t).$$ 
From the method of characteristics we have  $\ds \dot{x} = \frac{dx}{dt} =u(x,t)$, so that the modified transport equation \eqref{eq: transport A potential} can be written as
\begin{equation*}
\dot{B}=-f'(B)-B\, \frac{\partial u}{\partial x}.  
\end{equation*}
Using local analysis around $B=0$ and $B=1$ to respectively determine the appropriate ranges for $\gamma_1$ and $\gamma_2$, we conduct linear approximation to $u$ and estimate $\ds\frac{\partial u}{\partial x}$ as $a = a(x)$. This leads to the ODE given by 
\begin{equation}
\dot{B}=-f'(B)- aB.   \label{eq: simplified ODE}
\end{equation}

\noindent{\bf Determining a range for  $\gamma_1$}:
To determine an appropriate range for  $\gamma_1$, we use  local analysis around $B=0$ so that \eqref{eq: simplified ODE} becomes
\begin{equation}\label{eq: transport ODE1}
\dot{B}=-2\gamma_1 B -a B, \qquad B(0)=B_0<0,
\end{equation}
which has the analytical solution
\begin{equation}
\label{eq:Bsol}
B=B_0 e^{-(2\gamma_1+a)t}.
\end{equation}
Since we are considering the case where the numerical solution yields the out-of-range solution $A < 0$, we choose  the initial condition $B_0$ to be negative.  Our goal is to ``nudge'' the numerical solution so that $A$ moves back into range. Accordingly, we need  $B$ to be an increasing function, or equivalently $\dot{B} > 0$ in \eqref{eq: transport ODE1}.  Clearly this requires $\ds \gamma_1 > -\frac{a}{2}$. Observe from \eqref{eq:Bsol} that it is always true that $B < 1$ (in fact $B < 0$), so $A$ will never fall out of range near the value $1$ as long as  $\ds \gamma_1 > -\frac{a}{2}$.  

Imposing this constraint is straightforward. Since $A$ is initially within $[0,1]$ for the whole domain and $f'(A) = 0$ everywhere, we start by solving the non-modified transport equation \eqref{eq: transport A 1D}. Now suppose that at some later time there is a point in the domain for which the numerical scheme computes $A < 0$.  This is equivalent to $B_0 < 0$ in \eqref{eq: transport ODE1}, establishing the need to modify the transport equation by adding an extra forcing term $-f'(A)=-2\gamma_1A$ from \eqref{eq:forcing}. 

\begin{remark}
\label{rem:eulergamma1}
We could simplify the analysis by considering the one-step forward Euler approximation of \eqref{eq: transport ODE1},
\begin{equation}
\label{eq:gamma1upper}
B=-(2\gamma_1+a) B_0 \Delta t + B_0.    
\end{equation}
Based on the arguments above requiring $B$ to be an increasing function, we still need $\ds \gamma_1 > -\frac{a}{2}$. Using \eqref{eq:gamma1upper}, to ensure that $B \le 1$, we impose an upper bound for $\gamma_1$ as 
\begin{equation*}
\gamma_1 \leq -\frac{a }{2}-\frac{1-B_0}{2 B_0 \Delta t}, 
\end{equation*}
which yields the finite range for $\gamma_1$ given by
\begin{equation}
-\frac{a}{2} < \gamma_1 \leq -\frac{a }{2}-\frac{1-B_0}{2 B_0 \Delta t}.
\end{equation}
Note that the upper bound is not tight, because when $B_0$ is close to 0, which is in general the case, the term $\ds -\frac{1-B_0}{2 B_0 \Delta t}$ is a very large number. 
\end{remark}

\noindent{\bf Determining a range for  $\gamma_2$}:

Using a similar approach, we now consider the local analysis around $B=1$, corresponding to the case where $A$ goes out of range near the value $1$.  The ODE in \eqref{eq: simplified ODE} around $B=1$  reduces to 
\begin{equation}\label{eq: transport ODE2}
\dot{B}=-2\gamma_2(B-1) -a B, \qquad B(0)=B_0>1,
\end{equation}
which yields the analytical solution
\begin{equation}
\label{eq:gamma2analytic}
B=(B_0-\frac{2\gamma_2}{2\gamma_2+a}) e^{-(2\gamma_2+a)t}+\frac{2\gamma_2}{2\gamma_2+a}.
\end{equation} 
In this case, we seek $\gamma_2$ so that $B$ is decreasing, or equivalently $\dot{B} < 0$, which will ``nudge'' the numerical solution so that $A$ gets back in range of possible physical solutions, $[0,1]$.  Clearly, then, we require $\ds \gamma_2>-\frac{a B_0}{2(B_0-1)}$.  To ensure $B$ remains non-negative, so that we don't fall out of range on the left side of the solution interval, we first observe that \eqref{eq:gamma2analytic} can also be written as
\begin{equation}
B=(B_0-1) e^{-(2\gamma_2+a)t}+\frac{a}{2\gamma_2+a}(e^{-(2\gamma_2+a)t}-1)+1.
\end{equation}
It is immediately apparent that $B > 1$ for all $\gamma_2$ as long as $a\leq 0$. It is also possible to show that $B\geq 0$ for $a\leq 1$. In this regard we recall that $a$ is the approximation of $\ds \frac{\partial u}{\partial x}$, whose physical value is generally less than  $1$ m/s.   Thus we can conclude that $B \geq 0$ holds for all $\gamma_2$ for the given problem. 

\begin{remark}
\label{rem:eulergamma2}
As in \ref{rem:eulergamma1}, again we consider the one-step forward Euler approximation
\begin{equation}
B=-(2\gamma_2+a) B_0 \Delta t + 2\gamma_2 t +B_0.
\end{equation}
The requirements of decreasing $B$ with $B \ge 0$ lead to the range of $\gamma_2$ given by
\begin{equation}
\label{eq:rangegamma2}
-\frac{a B_0}{2(B_0-1)} < \gamma_2 \leq -\frac{a B_0}{2(B_0-1)}+\frac{B_0}{2(B_0-1) \Delta t}.
\end{equation}
\end{remark}

As before, we add the extra forcing term, in this case $-f'(A) = -2\gamma_2(A-1)$, whenever $A > 1$ results from the numerical solver. Finally, we also note that the same procedure is implemented for the ice thickness $h$ when the numerical solver causes it to become negative.

\section{Numerical experiments}
\label{sec:numerical}
We  provide three numerical experiments to illustrate the behavior of our proposed methods for the 1D sea ice simulation model. The first experiment is to corroborate the higher rate of convergence for WENO as compared to CD for smooth solutions. The capacity to resolve discontinuities (sharp features) in the  sea ice covers is verified in the second test, while the final example shows how the potential function is implemented in situations where the numerical solutions $A$ and $h$ fall out of range.

\subsection{Numerical convergence analysis}
\label{sec:convergence}

To assess the convergence rate of  WENO, we consider a test problem in a domain $\Omega=[0, \, 2000]$ km with a known analytical solution.\footnote{This is, of course, generally not the case as the sea ice model has no known solutions.} We construct our test cases by adding appropriate extra forcing terms to the governing equations. Specifically, we introduce to the right hand side of each equation of \eqref{eq:1DVPmodel} forcing terms which are obtained by plugging into the left hand side terms corresponding to the following proposed solutions: 
\begin{eqnarray}
u_{\text{true}} &=&(\sin(2\pi x/(2\times 10^{6}) + 5t/518400-\pi/2)+1)\times 0.001+0.2, \nonumber\\
h_{\text{true}} &=&(\sin(2\pi x/(2\times 10^{6}) + 5t/518400-\pi/2)+1)+0.1,\nonumber\\
A_{\text{true}} &=&(\sin(2\pi x/(2\times 10^{6}) + 5t/518400-\pi/2)+1)\times 0.15+0.7.
\label{eq:testsolution}
\end{eqnarray}
Observe that the values in \eqref{eq:testsolution} of $u_{\text{true}}$, $h_{\text{true}}$ and $A_{\text{true}}$ are consistent in magnitude to their corresponding true physical values. Also observe that while we construct extra forcing terms from the left hand side of \eqref{eq: momentum u 1D}, the wind forcing is canceled out and therefore does not affect the convergence results.  Initial conditions for the system are obtained by plugging $t = 0$ into \eqref{eq:testsolution}.

The total simulation time is $T=5$ s.  The time step $\Delta t=10^{-4}$ s is intentionally chosen to be small enough to ensure that the time discretization error does not affect the convergence rates. It furthermore allows us to conduct convergence tests directly on the VP sea ice model for both the WENO and CD explicit spatial discretization schemes  without the usual concern for the stability issue associated with explicit methods.

Table \ref{table: convergence} compares the relative $\ell_2$ errors for increasing resolutions with each spatial discretization choice. We observe second-order convergence for all three variables for the CD case, which is consistent with the standard convergence analysis results for CD schemes.  By employing the WENO scheme in the ideal case, one would expect to obtain sixth-order convergence for the velocity $u$ and fifth-order convergence for both the ice thickness $h$ and ice concentration $A$. However, due to the complexity and non-linearity of the sea ice model, coupled with the fact that the added extra forcing terms are not being updated in the stages of TVRK3 time integration, theoretical accuracy is unlikely to be obtained. We still observe higher-order convergence for all three variables as compared to the CD results. In addition, a direct comparison of the error magnitudes for all three variables indicates that the WENO scheme indeed provides more accurate results than those obtained using CD, noting that WENO appears to be mainly affected by round-off error at $10$ km resolution.

\begin{table}[ht]
\begin{center}
\begin{tabular}{  c | c | c | c  | c | c  | c }
\hline
\multicolumn{7}{c}{CD TVRK3} \\
\hline
resolution $\Delta x$ & u error & u rate & h error & h rate & A error & A rate\\ \hline
40 km& 2.6655e-06 &        & 4.4967e-09 &         & 1.0362e-09 & \\ \hline
20 km& 6.6698e-07 & 1.9987 & 1.1247e-09 & 1.9992 & 2.5920e-10 & 1.9991\\ \hline
10 km& 1.6692e-07 & 1.9984 & 2.8120e-10 & 1.9998 & 6.4883e-11 & 1.9981\\ \hline
\end{tabular}

\medskip
\begin{tabular}{  c | c | c | c  | c | c  | c }
\hline
\multicolumn{7}{c}{WENO TVRK3} \\
\hline
resolution $\Delta x$ & u error & u rate & h error & h rate & A error & A rate\\ \hline
40 km& 5.2407e-07 &        & 1.3483e-11 &        & 8.8200e-12 & \\ \hline
20 km& 2.1769e-08 & 4.5894 & 5.8573e-13 & 4.5248 & 9.2062e-13 & 3.2601\\ \hline
10 km& 8.3211e-10 & 4.7093 & 8.8497e-14 & 2.7265 & 5.5688e-13 & 0.7252 \\ \hline
\end{tabular}
\end{center}
\caption{
A comparison of CD and WENO spatial discretization errors for increasing resolution.}\label{table: convergence}
\end{table}

\subsection{A simulation of sea ice with sharp features}
\label{sec:sharpfeatures}

In this example we test the performance of the WENO scheme on a simulation of a sea ice cover with sharp features.  To better capture the solution behavior near the discontinuity region while maintaining periodic boundary conditions, the structure of ice is designed such that relatively solid ice covers both ends of the domain and a very thin layer of ice is in the center of the domain. This is realized via a discontinuous setting on the initial conditions of ice thickness and ice concentration given by
\begin{equation}
\begin{gathered}
\ds u=0 \text{ m/s} \qon [0, 2000] \text{ km},\\
\ds h= \begin{cases}
0.01 \text{ m} \qon [400, 1600] \text{ km},\\
2 \text{ m} \qon [0,400] \cup [1600, 2000] \text{ km},
\end{cases}\\
\ds A= \begin{cases}
0 \qon [400, 1600] \text{ km},\\
\ds 0.8 \qon [0,400] \cup [1600, 2000] \text{ km}.
\end{cases}
\label{eq:initialconditions}
\end{gathered}
\end{equation}
For the external forcing in \eqref{eq: tau a} we impose uniform constant wind forcing $u_a=10$ m/s. 


As the WENO and CD schemes yield theoretically different convergence rates, for a direct comparison, we also consider the linear WENO scheme, for which the nonlinear weights are replaced by linear ones of the same order accuracy. Note that this is equivalent to using an upstream centered scheme (upstream in time, centered in space), \cite{Jiang96}.  Due to the combined stencil, the highest possible order of accuracy is obtained in smooth regions. The results are oscillatory near discontinuities, however. 


\subsubsection{Simulation results on the VP model}
\label{sec:simulateVP}
We first discuss the results for the VP model in \eqref{eq:1DVPmodel}. The simulation is run with a spatial resolution of $\Delta x=10$ km and time step $\Delta t=1$ s for a total simulation time of $1$ hour ($3600$ s).

\begin{figure}[ht]
\centering
\includegraphics[width=0.32\textwidth]{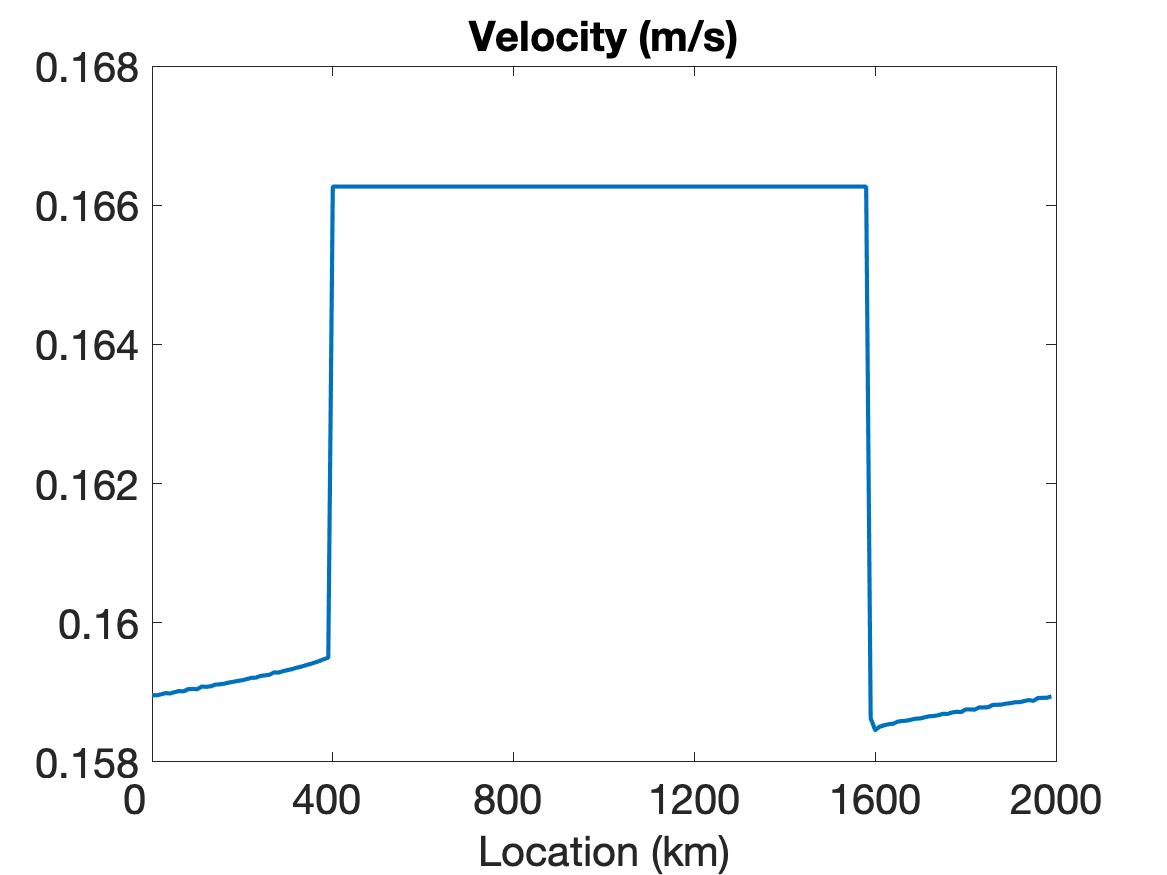}
\includegraphics[width=0.32\textwidth]{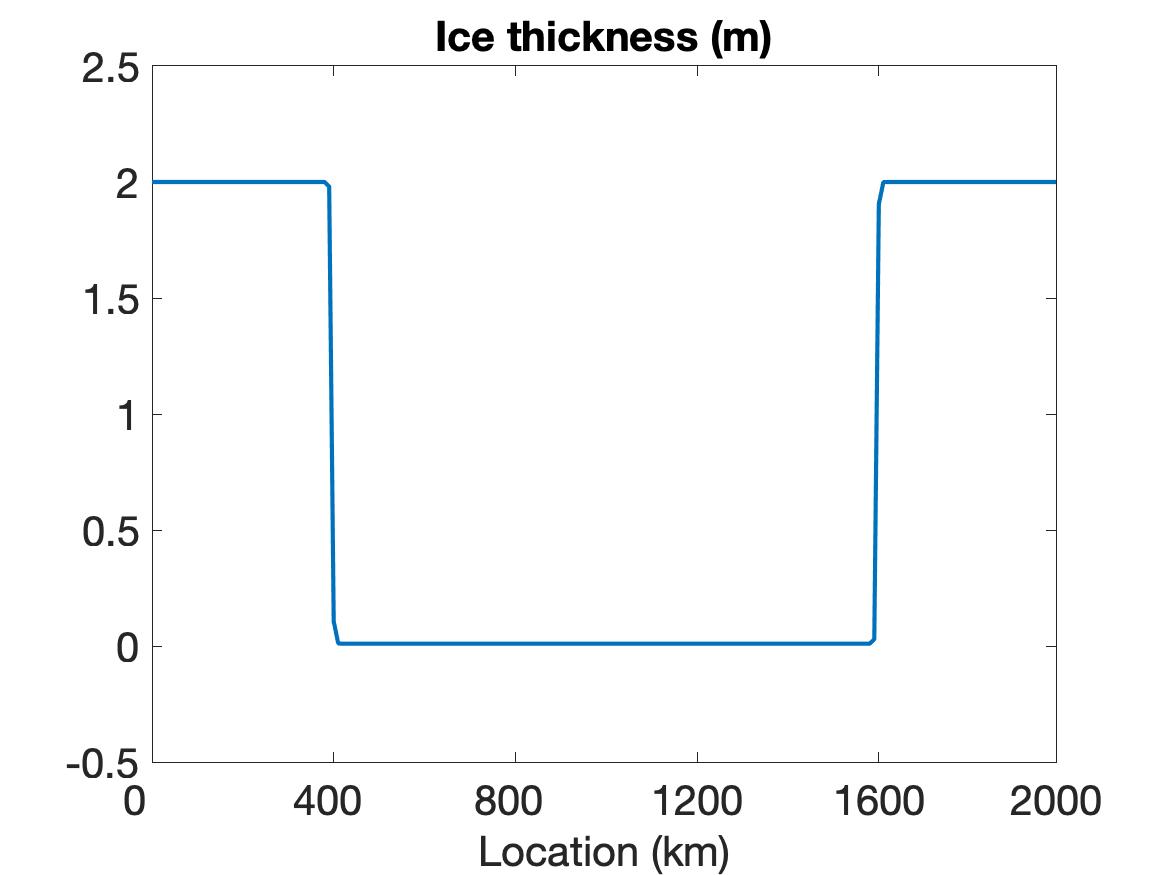}
\includegraphics[width=0.32\textwidth]{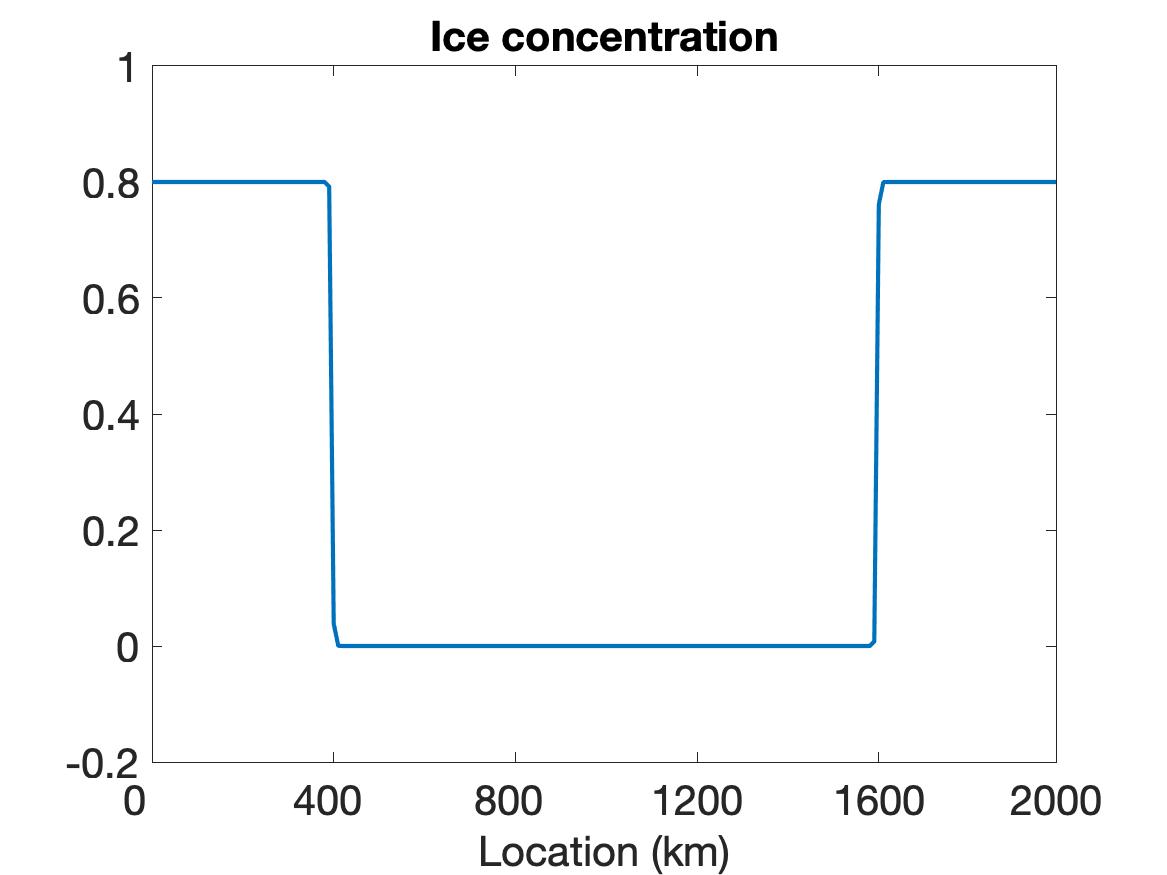}\\[1ex]
\includegraphics[width=0.32\textwidth]{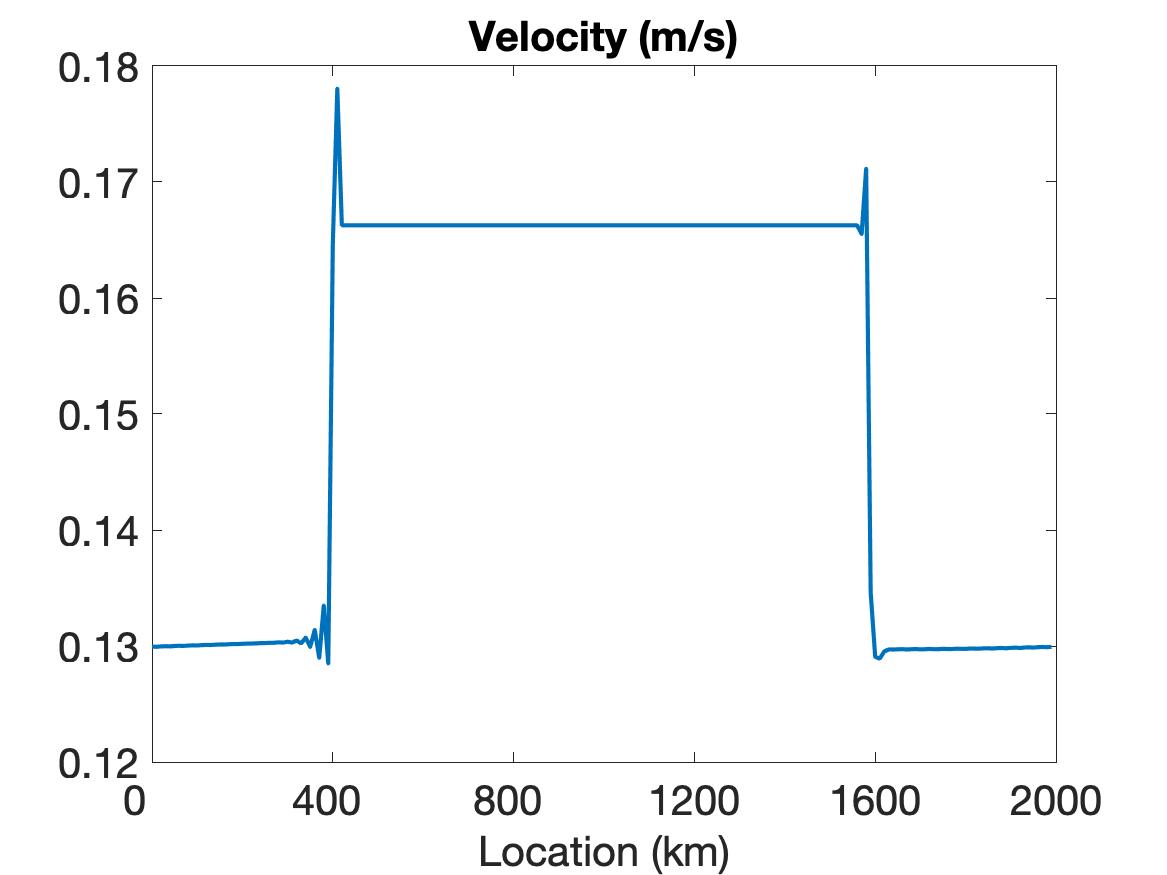}
\includegraphics[width=0.32\textwidth]{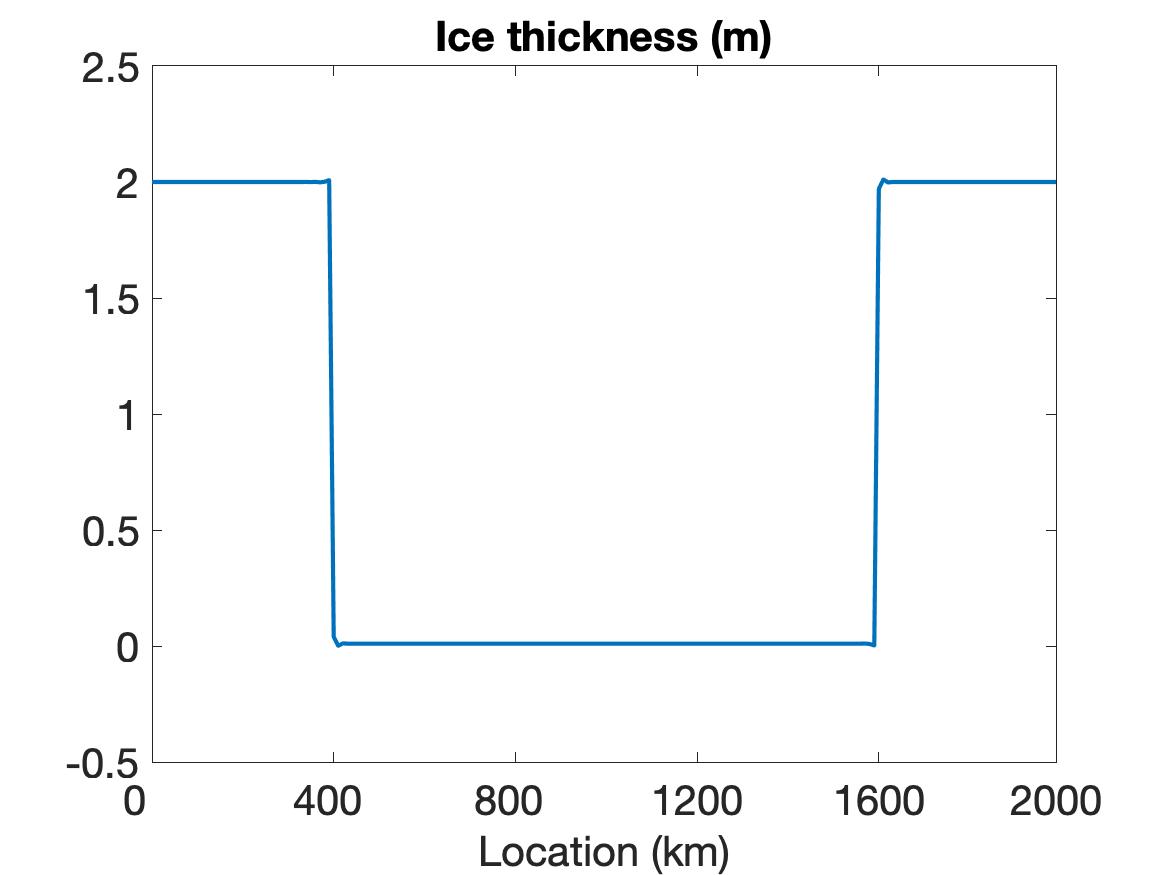}
\includegraphics[width=0.32\textwidth]{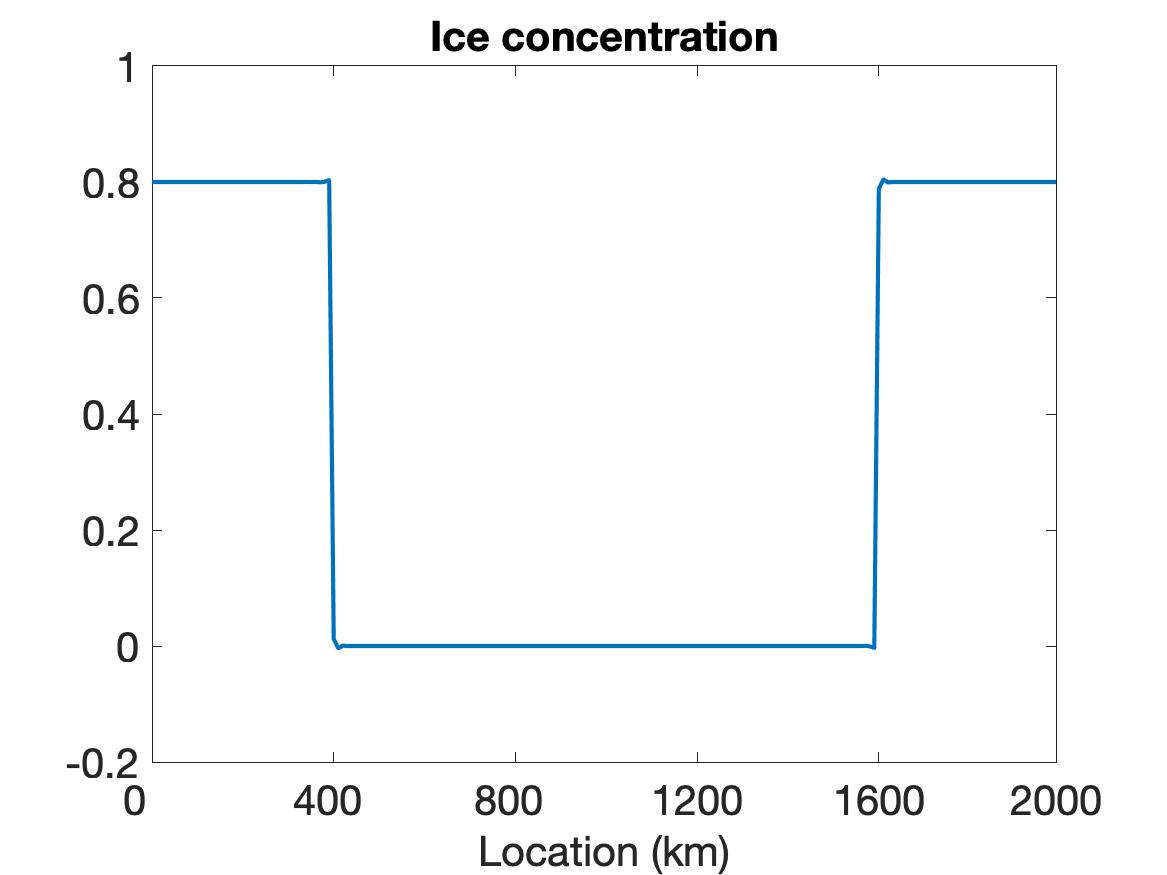}\\[1ex]
\includegraphics[width=0.32\textwidth]{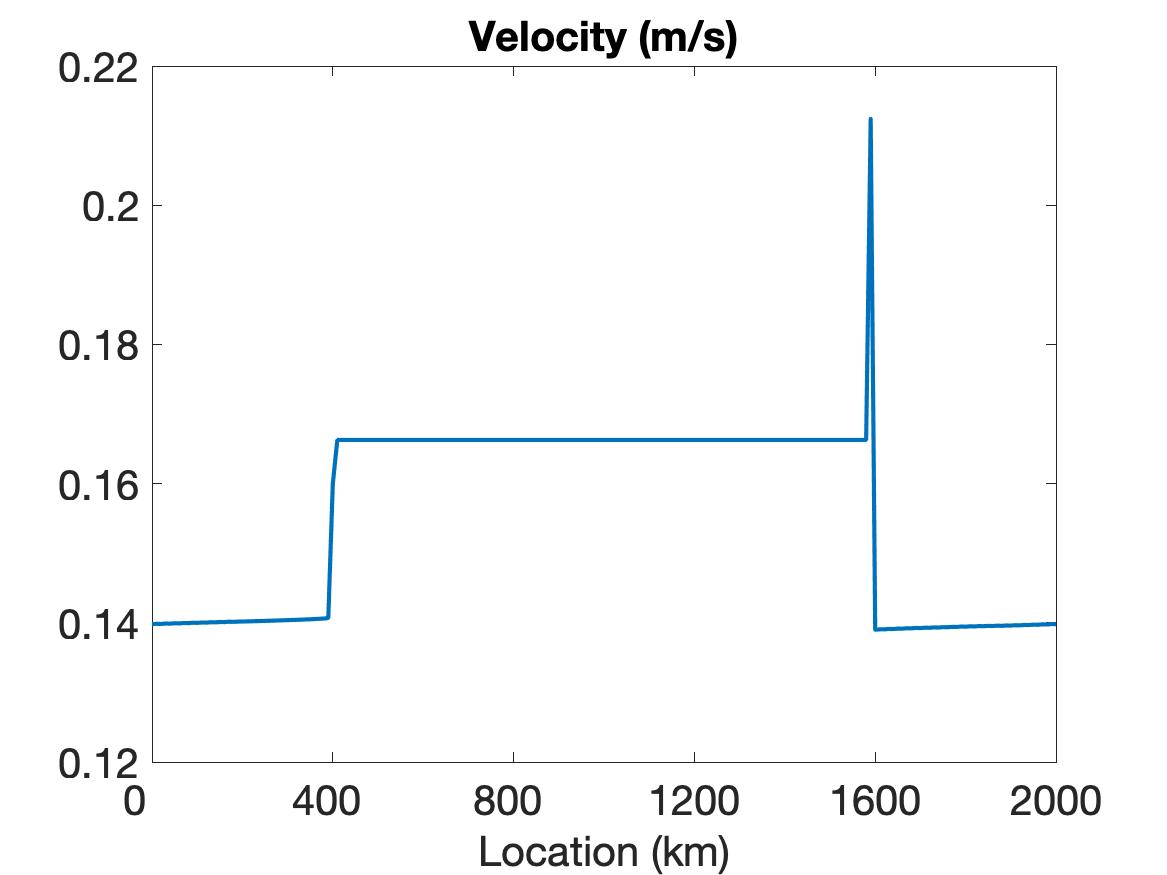}
\includegraphics[width=0.32\textwidth]{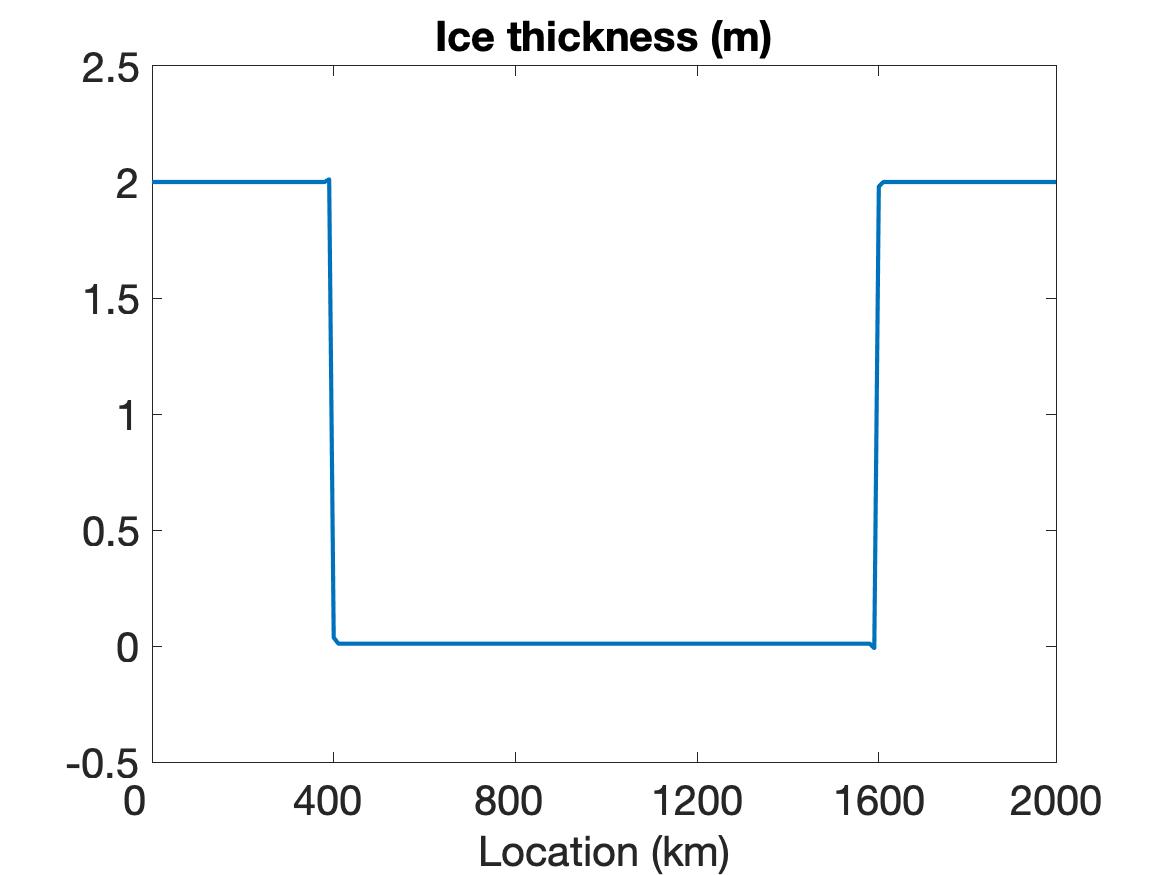}
\includegraphics[width=0.32\textwidth]{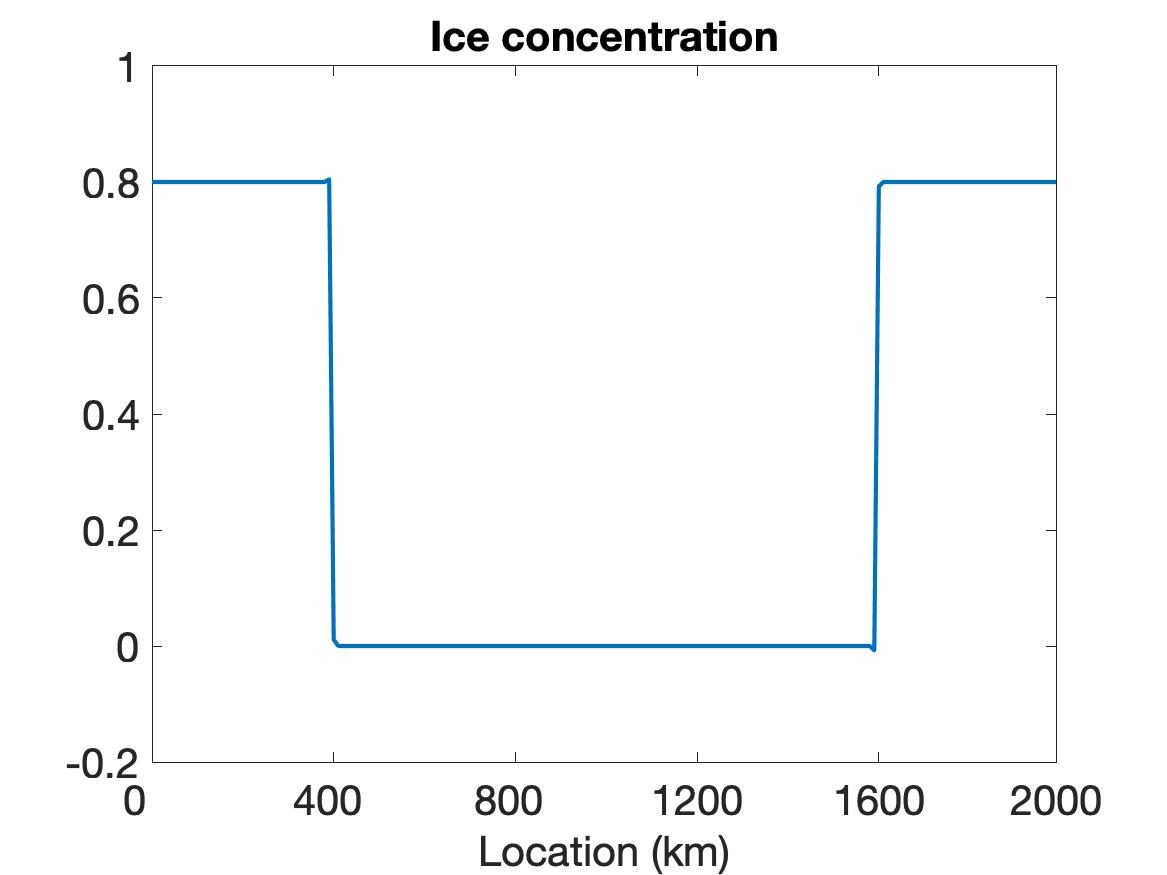}
\caption{Simulation of sea ice with sharp features on VP model. (Top) solution plots using WENO at $1$ hour; (middle row) solution plots using linear WENO at $2000$ s; (bottom) solution plots using CD at $2332$ s. (Left) velocity $u$; (middle column) ice thickness $h$; (right) ice concentration $A$.}
\label{fig: VP discontinuous}
\end{figure}

Figure \ref{fig: VP discontinuous} compares the results using WENO (top row), linear WENO (middle row), and CD (bottom row) spatial discretizations for the simulation of the sea ice model with sharp features as constructed using the initial conditions given in \eqref{eq:initialconditions}. Observe that the solutions for each variable  $u$ (left column), $h$ (middle column) and $A$ (right column) are discontinuous. The solution plots demonstrate that only WENO maintains a sharp non-oscillatory solution for the velocity $u$ in each sub-region, with a sharp overshoot occurring in the CD velocity profile.  We also note that while the WENO solution is plotted at the final time of $1$ hour, the solutions for the linear WENO and CD are presented at $2000$ s and at $2332$ s, respectively, since the oscillations eventually cause these solutions to blow up. There is less distinction between the methods in the ice thickness and ice concentration solutions, which all retain the initial profiles while  moving slightly toward the ends of the domain due to the exerted wind forcing.  However, since they are coupled with velocity, the CD and linear WENO solutions will also blow up before the final time. 

\subsubsection{Simulation results on the EVP model}
\label{sec:simulateEVP}

For the EVP solver described in Section \ref{sec:EVP}, we run the simulation with spatial resolution $\Delta x=10$ km and time step $\Delta t=10$ s with $1000$ sub-cycling steps for a total simulation time of $1$ hour.

\begin{figure}[ht]
\centering
\includegraphics[width=0.32\textwidth]{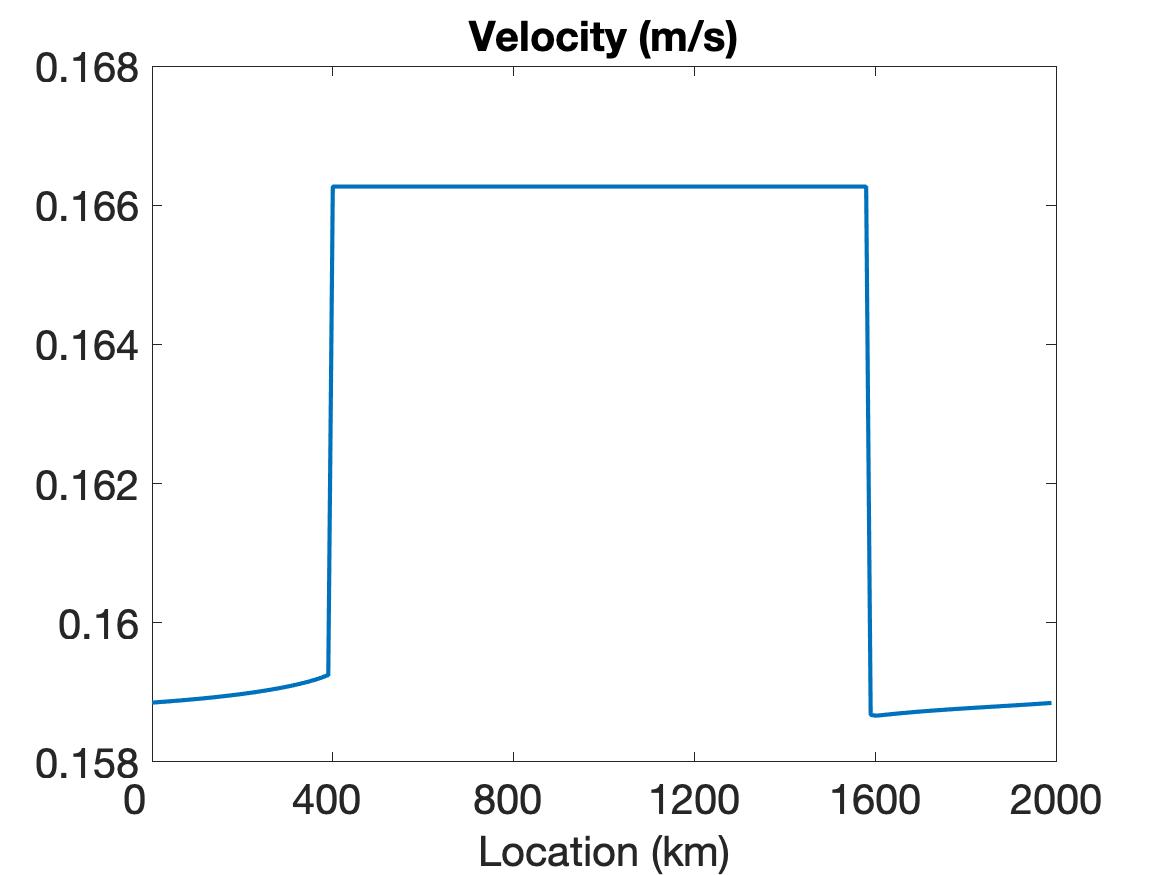}
\includegraphics[width=0.32\textwidth]{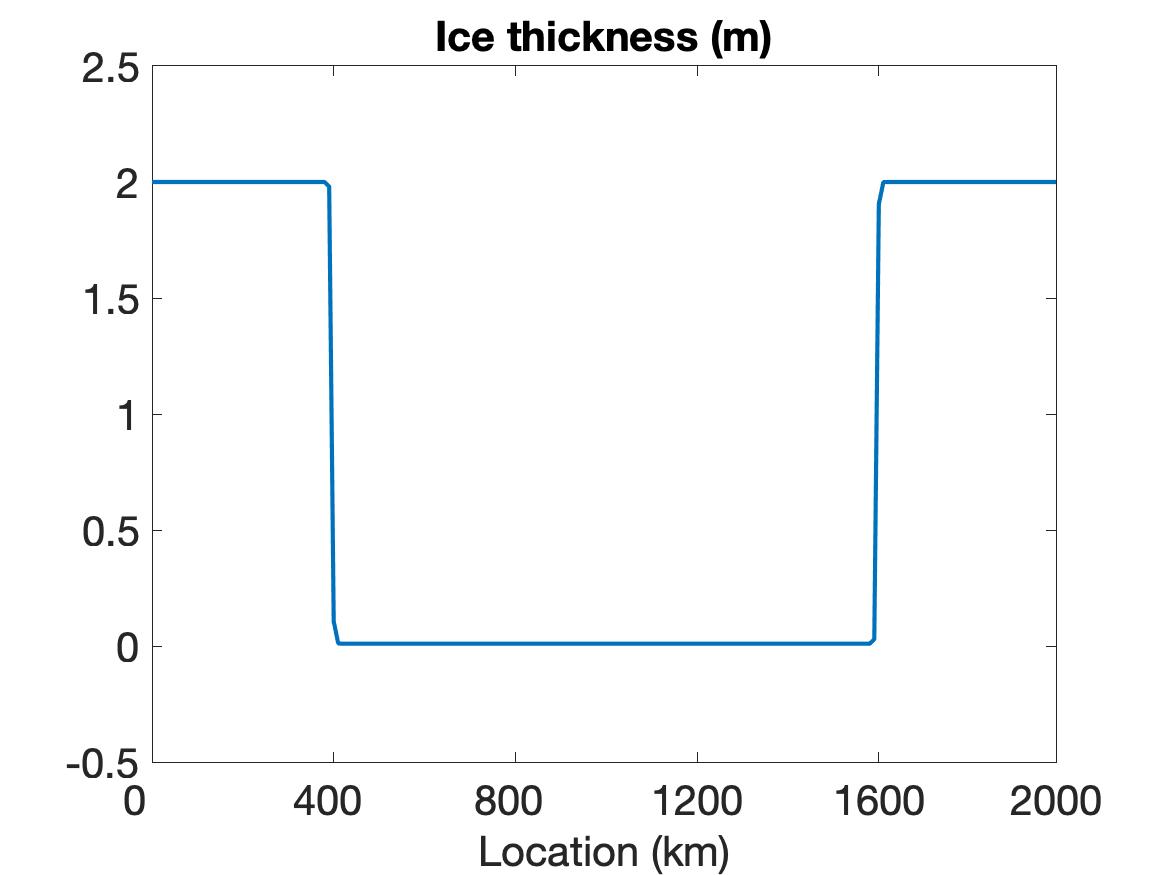}
\includegraphics[width=0.32\textwidth]{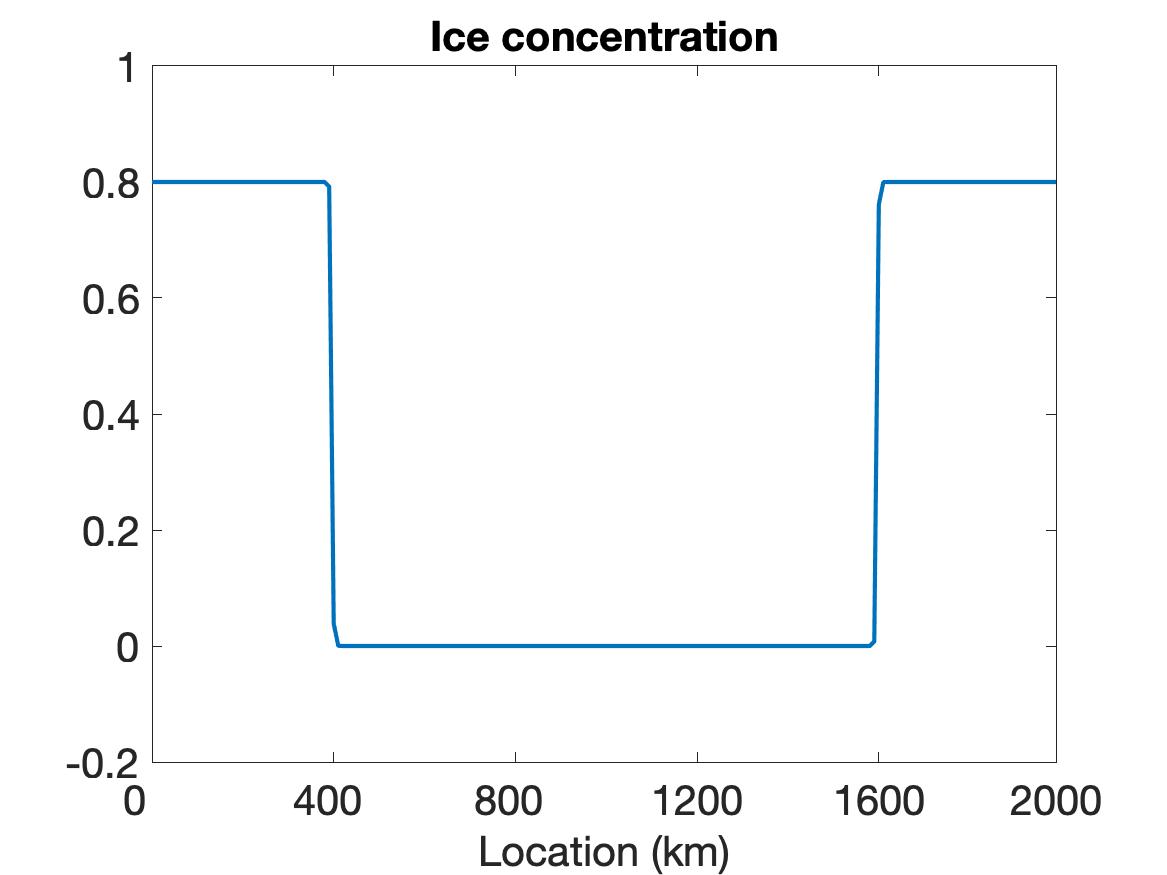} \\[1ex]
\includegraphics[width=0.32\textwidth]{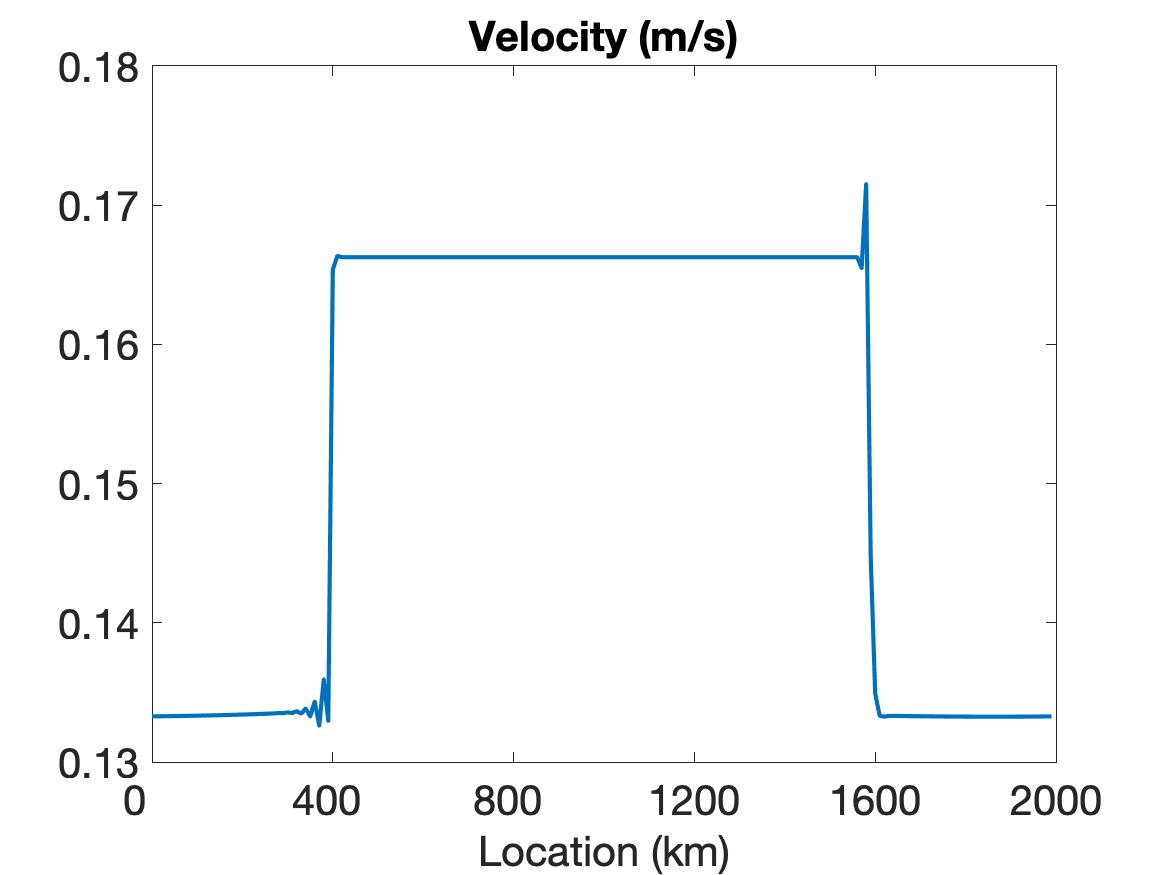}
\includegraphics[width=0.32\textwidth]{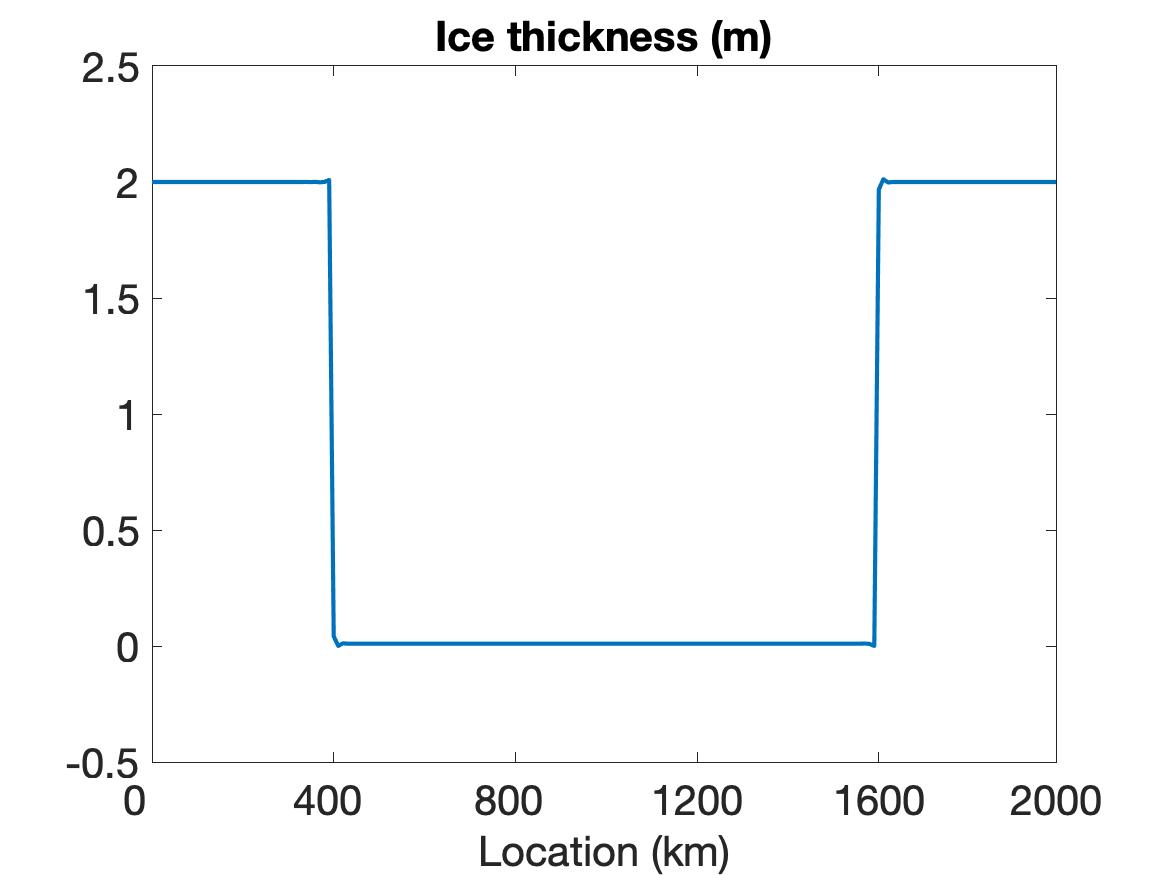}
\includegraphics[width=0.32\textwidth]{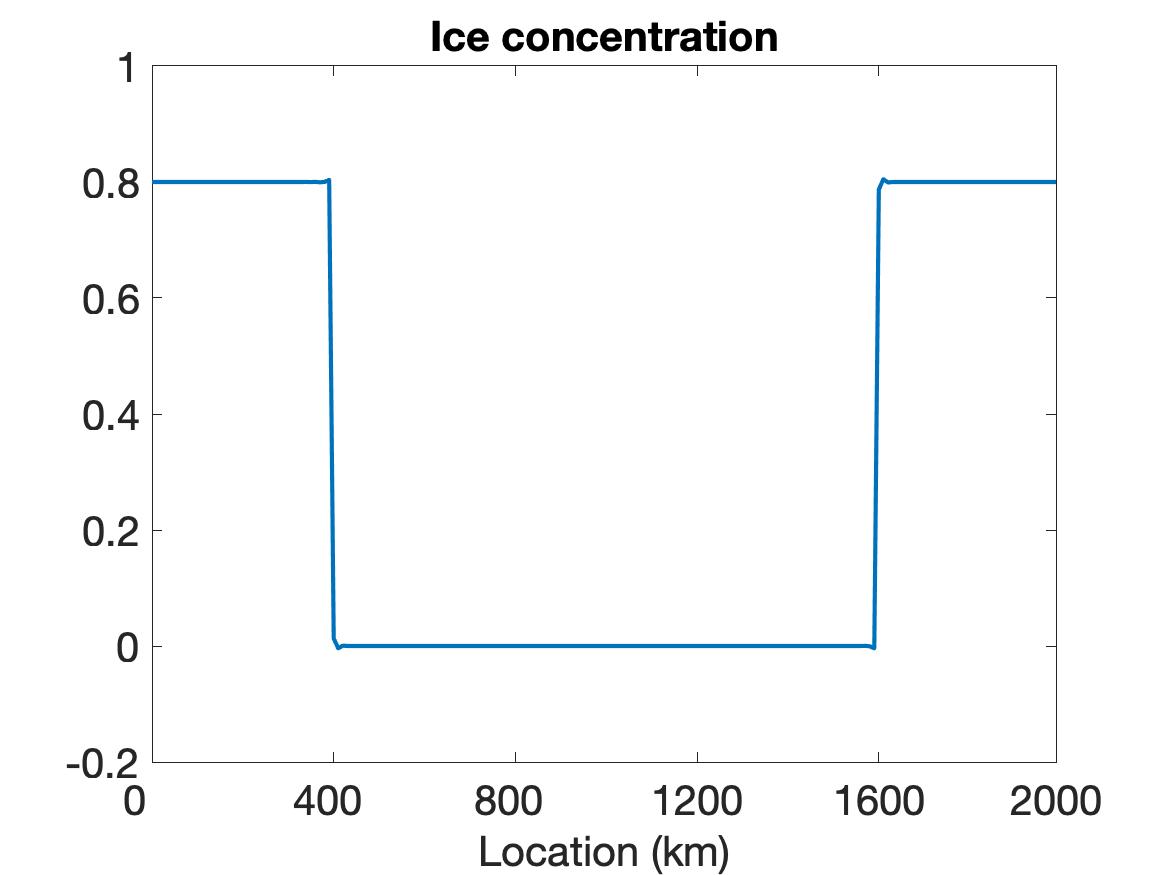}\\[1ex]
\includegraphics[width=0.32\textwidth]{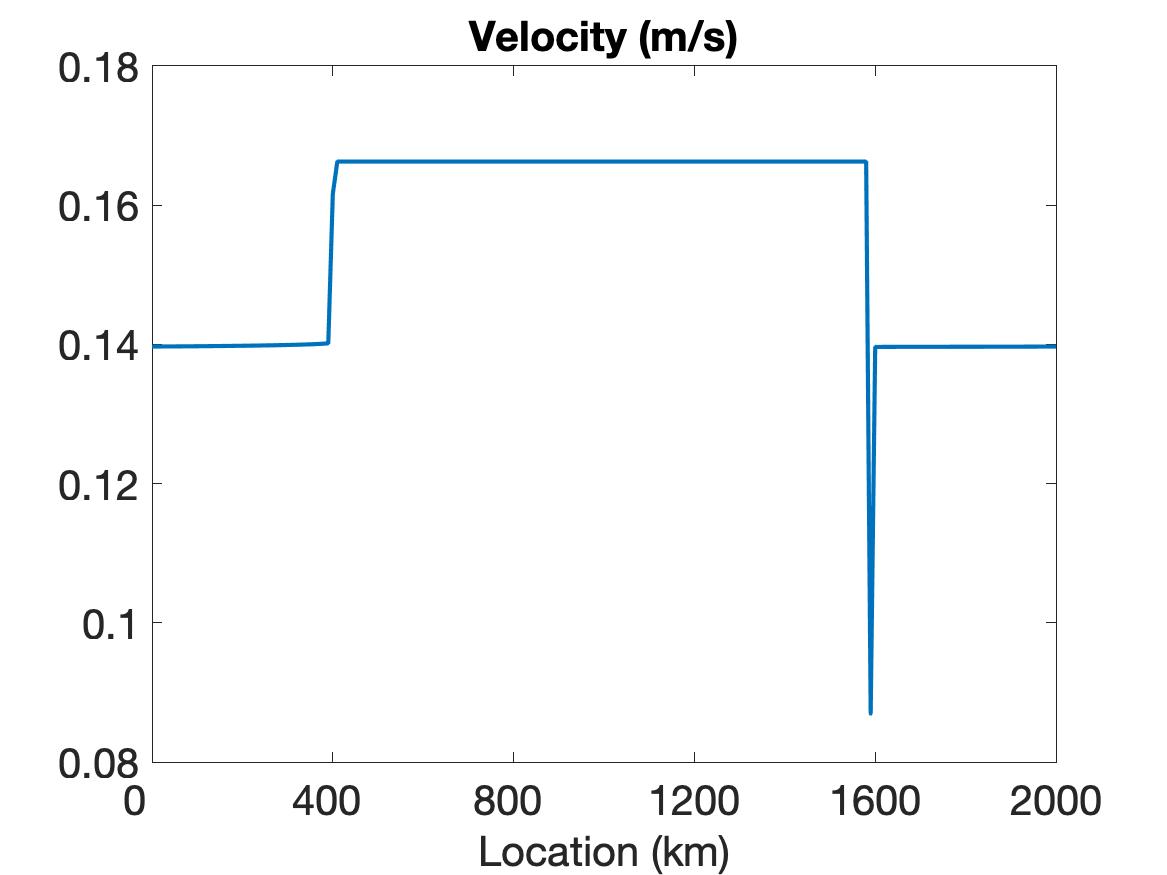}
\includegraphics[width=0.32\textwidth]{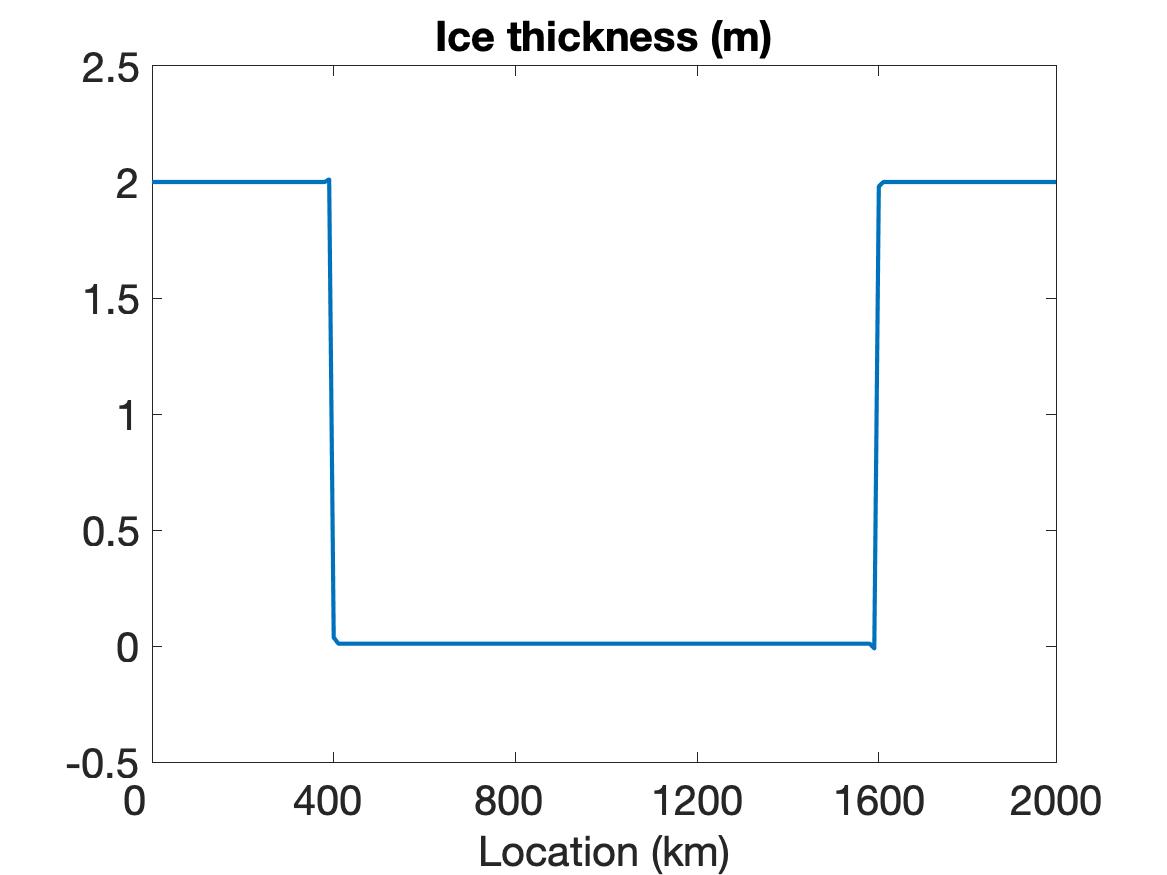}
\includegraphics[width=0.32\textwidth]{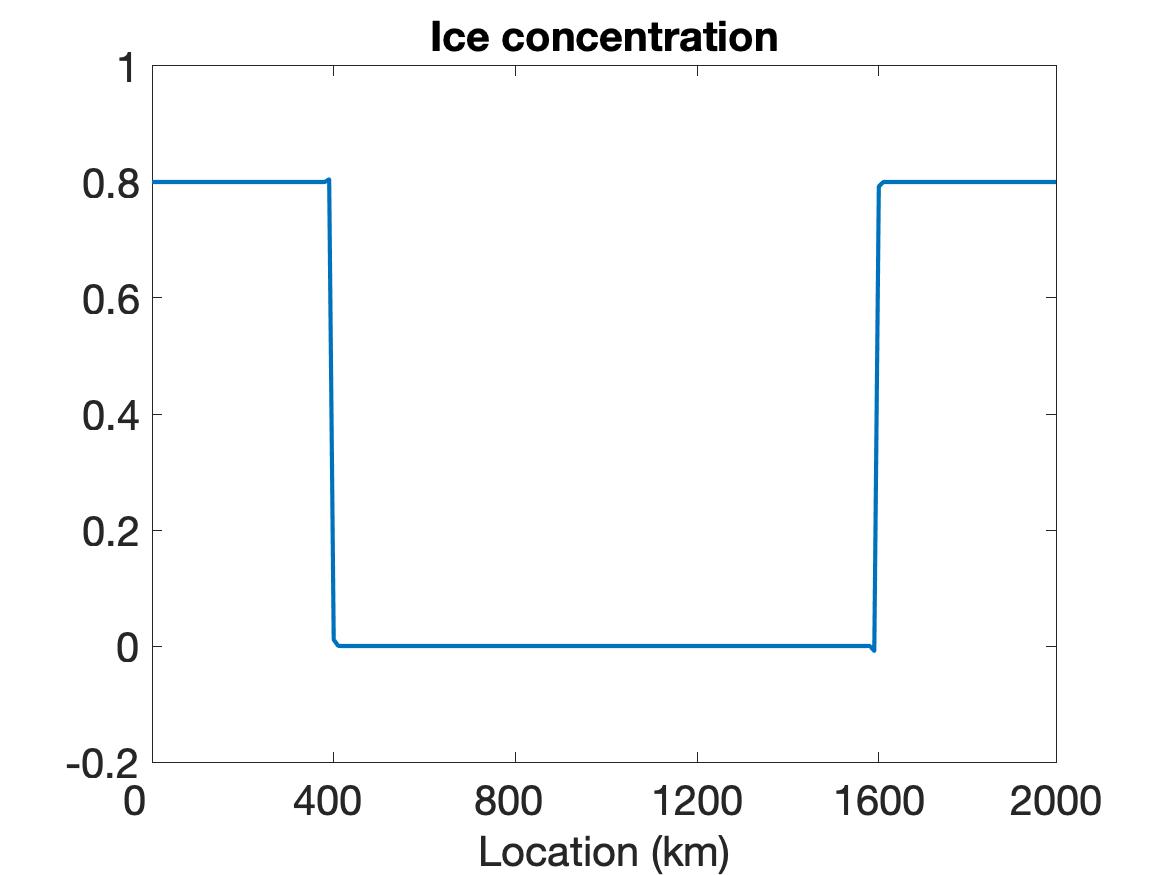}
\caption{Simulation of sea ice with sharp features on EVP model. (Top) solution plots of $u$, $h$ and $A$ using WENO at $1$ hour; (middle row) solution plots of $u$, $h$ and $A$ using linear WENO at $2110$ s; (bottom) solution plots of $u$ (left) during sub-cycling between $2340$ s and $2350$ s; $h$ (middle) at $2340$ s; $A$ (right) at $2340$ s using CD .} 
\label{fig: EVP discontinuous}
\end{figure}

As in Figure \ref{fig: VP discontinuous}, Figure \ref{fig: EVP discontinuous} displays the solution for the initial conditions given in \eqref{eq:initialconditions}  obtained by WENO (top row), linear WENO (middle row) and CD (bottom row) spatial discretizations. Once again, the simulation can only reach the final time of $1$ hour using the WENO scheme.  Observe that the results for the EVP and VP models are nearly identical, with no oscillatory behavior near discontinuities. The linear WENO solution is shown at time $2110$ s, where again we see oscillations in the velocity profile. The bottom row (left) shows the velocity during the sub-cycling iteration between $2340$ s and $2350$ s. The ice thickness and ice concentration at $2340$ s are shown in the bottom-middle and bottom-right, respectively. We present the velocity profile during the sub-cycling stage to capture the undershoot that occurs within the sub-cycling  -- it is not detectable outside the sub-cycling for this case. This undershoot eventually leads to the solution blowing up before reaching the final time.

The simulation results for both the VP and EVP models lead us to conclude that while we are able to properly resolve the discontinuities and obtain a stable solution using WENO, this cannot be accomplished using either the CD or the traditional higher-order (linear WENO) schemes.

\subsection{Incorporating the potential function method into the solver}
\label{sec:numtestpotential}
We now test the model for which no exact solution is known. The main goal of this numerical test is to determine how out-of-range issues, namely $A < 0$, $A > 1$, or $h < 0$,  may be effectively mitigated using the potential method described in Section \ref{sec:potential}.  

We use the backward Euler time integration with CD spatial discretization and JFNK to solve the momentum equation, along with TVRK3 with CD spatial discretization scheme for transport equations.
We note that in choosing to use an  implicit time-stepping method for solving \eqref{eq: momentum u 1D}  we avoid issues concerning stability. Moreover, with regard to the transport problem, we note that the WENO scheme does not yield out-of-range negative solutions for either $A$ or $h$. Hence to determine the efficacy of the potential method we apply the CD scheme (using explicit time-stepping) to the transport equations \eqref{eq: transport h 1D} and \eqref{eq: transport A 1D}.  As discussed in \ref{sec:WENO}, we also use the one-dimensional version of the staggered Arakawa C-grid \cite{Arakawa77}.
 

We run all the experiments on a domain of length $2000$ km for a total integration time of 6 days. We choose a spatial resolution of $20$ km  and $90$ s as the time step. The initial conditions are constant throughout $[0,2000]$ km and given by
\begin{equation*}
u(x,0) =0 \ \text{m/s}, \quad h(x,0)= 1 \ \text{m}, \quad A(x,0) =0.9. 
\end{equation*}
We impose Dirichlet boundary conditions so that 
\begin{equation*}
u(0,t) = u(2000,t) = 0 \text{ m/s}.
\end{equation*}

Observe that since the ocean is considered to be at rest in this model, the only variable external forcing term is the wind, which is imposed uniformly as a constant given by
\begin{equation*}
u_a(x,t) =10  \ \text{m/s}. 
\end{equation*}

For imposing the boundary conditions numerically, we have $u_0 = 0$ (similarly $u_N = 0$) so that
$$ (uh)_0=u_0 h_0=0,$$
leading to
$$\{d(uh)\}_{\frac{1}{2}}=\frac{(uh)_{1}-(uh)_0}{\Delta x}=\frac{(uh)_{1}}{\Delta x}$$
for the transport equation \eqref{eq: transport h 1D} 
(similarly for $\{d(uh)\}_{N-\frac{1}{2}}$).  The Dirichlet boundary conditions also yield analogous equations for $A$ at the boundaries  in \eqref{eq: transport A 1D}.

For comparison purposes we first run the simulation without applying the potential function method. Throughout the 6-day simulation, we find that the largest $A$ value is 1.0540, the smallest $A$ value is -0.1445, and the smallest $h$ value is -0.1606, which are all out of range.  

To employ the potential function method, we first  simulate the model (without incorporating the potential function method into the transport equations) until the time {$T_1$} for which  $\ds \min_x\{A(x,{T_1})\} < 0$.\footnote{We only do this process one time for each out-of-range situation, $A < 0$, $A > 1$, and $h < 0$.  The potential function parameters then remain fixed for all time.} The potential function variables corresponding to \eqref{eq: transport ODE1} are then determined as 
$$ a \approx \frac{\partial u}{\partial x}, \quad B_0=A.$$
While $B_0$ is determined directly from the numerical implementation of the scheme, as previously mentioned after \eqref{eq: simplified ODE}, $a$ is computed as a linear approximation of $\ds \frac{\partial u}{\partial x}$. 

Following the discussion in Section \ref{sec:potential}, we then find a uniform bound for $\gamma_1$ by computing 
$$\gamma_{1,\min} = \ds \max_x \Big\{-\frac{a}{2}\Big\}, \quad\quad \gamma_{1,\max} = \ds \min_x \Big\{-\frac{a }{2}-\frac{1-B_0}{2 B_0 \Delta t}\Big\}.$$ 
The process is similar for determining $\gamma_2$ for when $\ds \max_x\{A(x,{T_2})\} > 1$, ${T_2} > 0$,  and in our experiment we obtain
$$\ds \gamma_1 \in (3.6682 \times 10^{-7}, \ 786.4101), \quad\quad \gamma_2 \in (0.0075, \ 273.1214).$$  
Finally the corresponding range for the parameter $\gamma$ associated with $h < 0$ is 
$$\gamma \in (3.6682 \times 10^{-7}, \ 707.7696).$$ 
Based on these results and the  discussion in Section \ref{sec:potential}, we choose
$$ \gamma_1=10^{-3}, \qquad \gamma_2=10^{-2}, \qan \gamma=10^{-3}.$$
We then run the remainder of the simulation, up until final time $T = 6$ days, with the potential function method now incorporated into the transport equations.

\begin{figure}[ht]
\begin{center}
\includegraphics[width=0.32\textwidth]{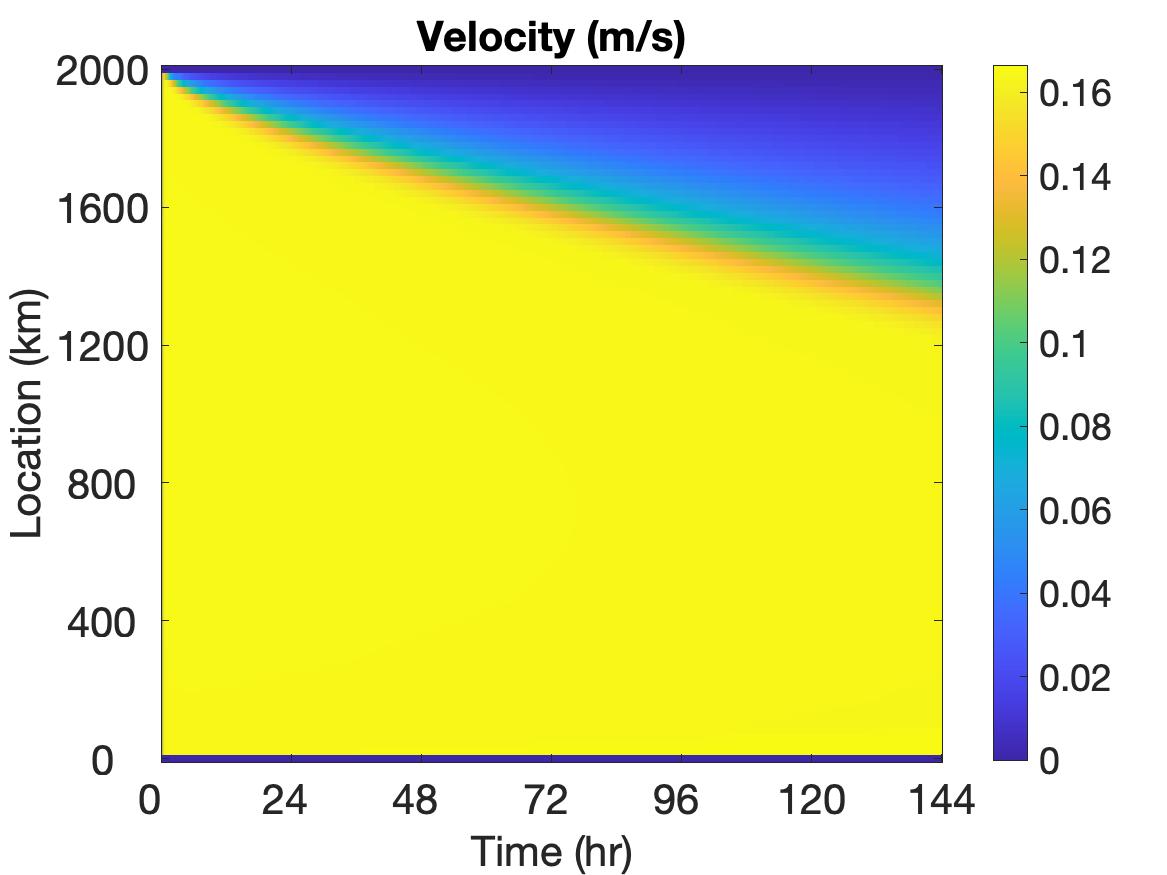}
\includegraphics[width=0.32\textwidth]{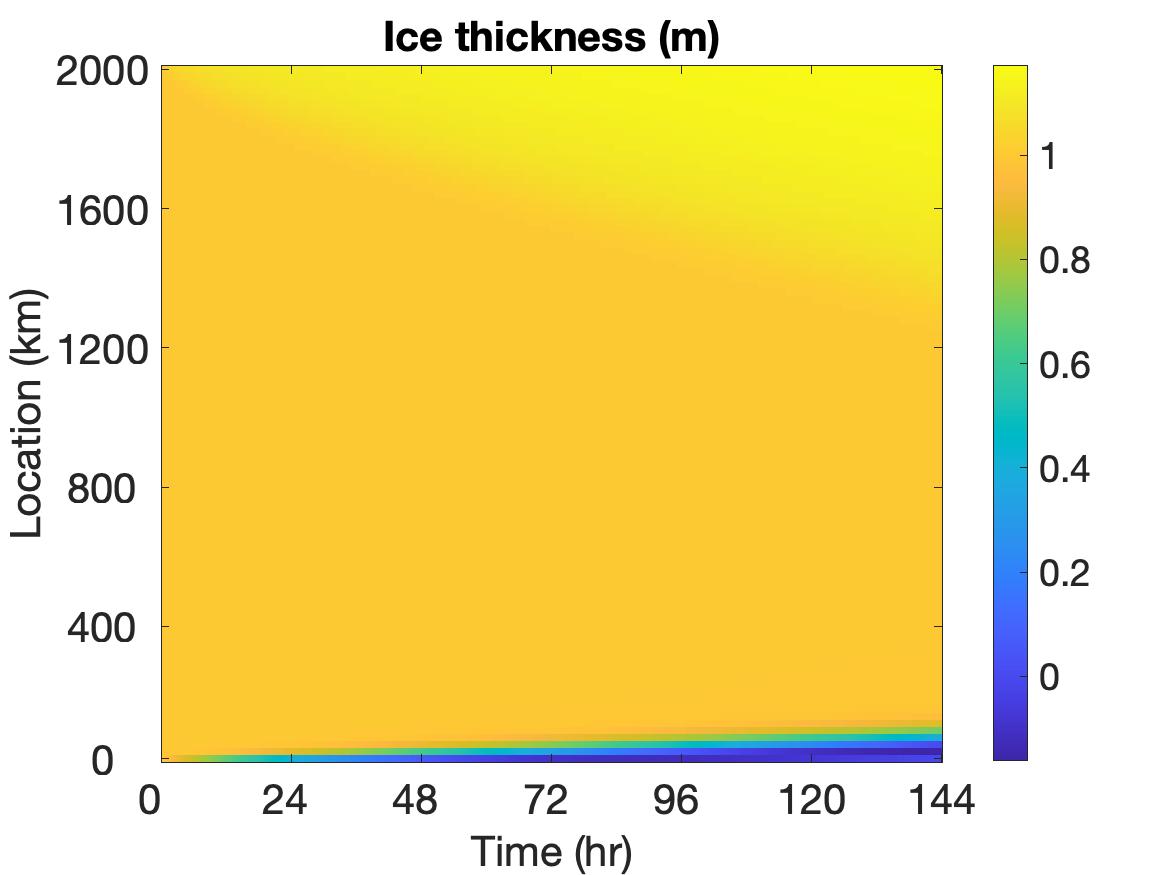}
\includegraphics[width=0.32\textwidth]{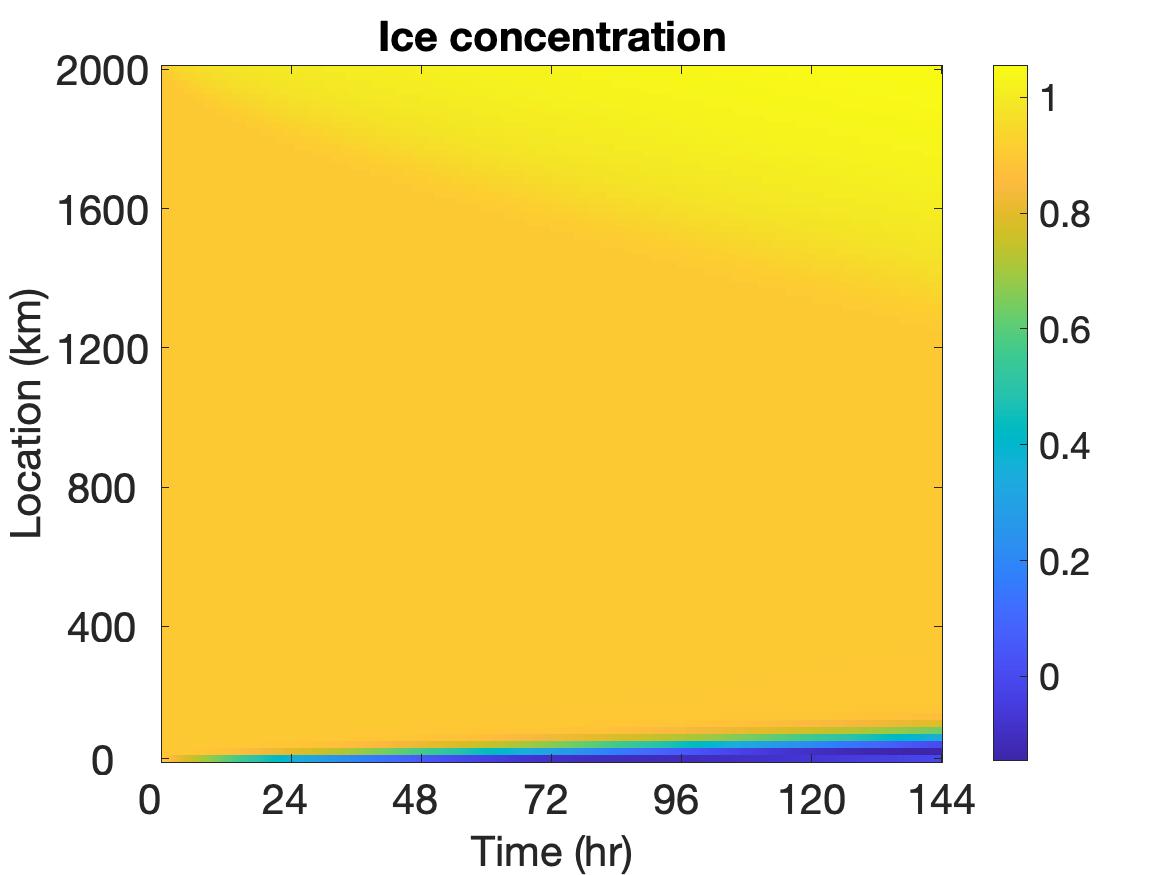}\\[1ex]
\includegraphics[width=0.32\textwidth]{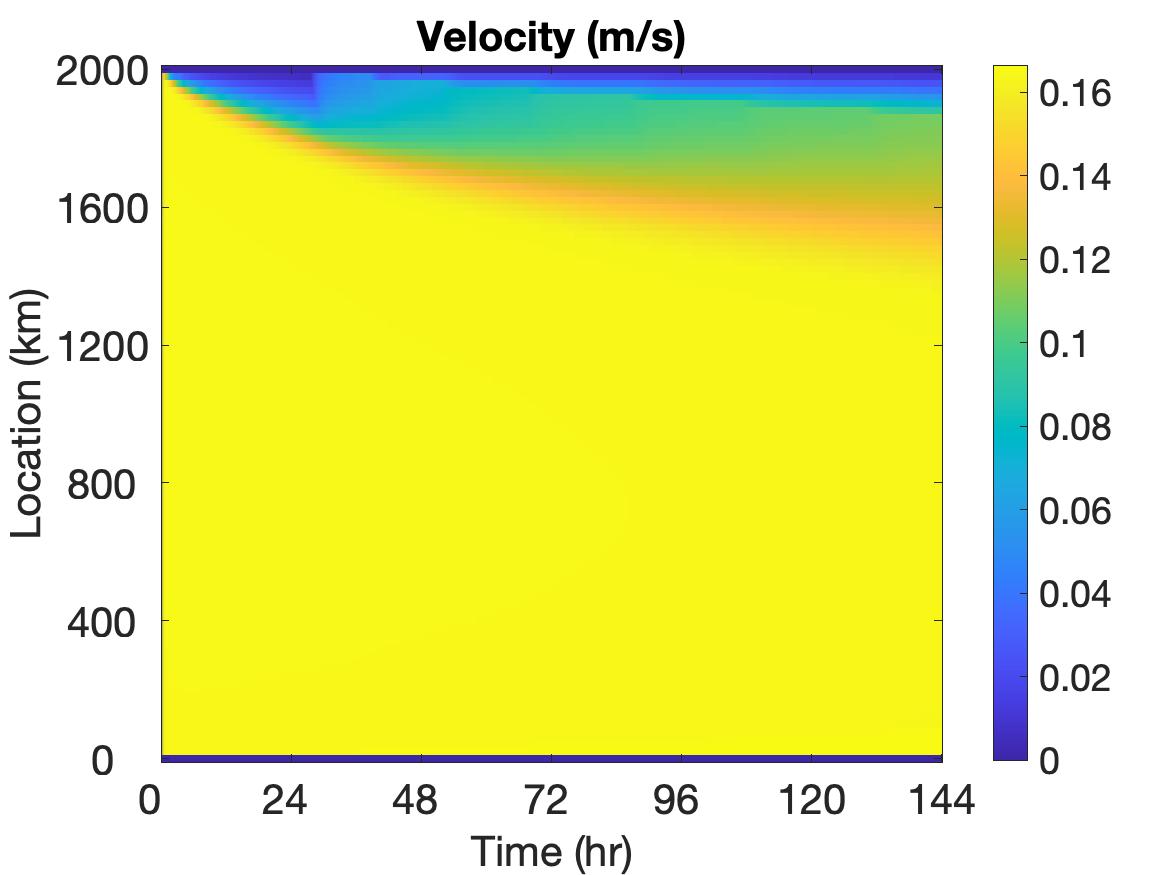}
\includegraphics[width=0.32\textwidth]{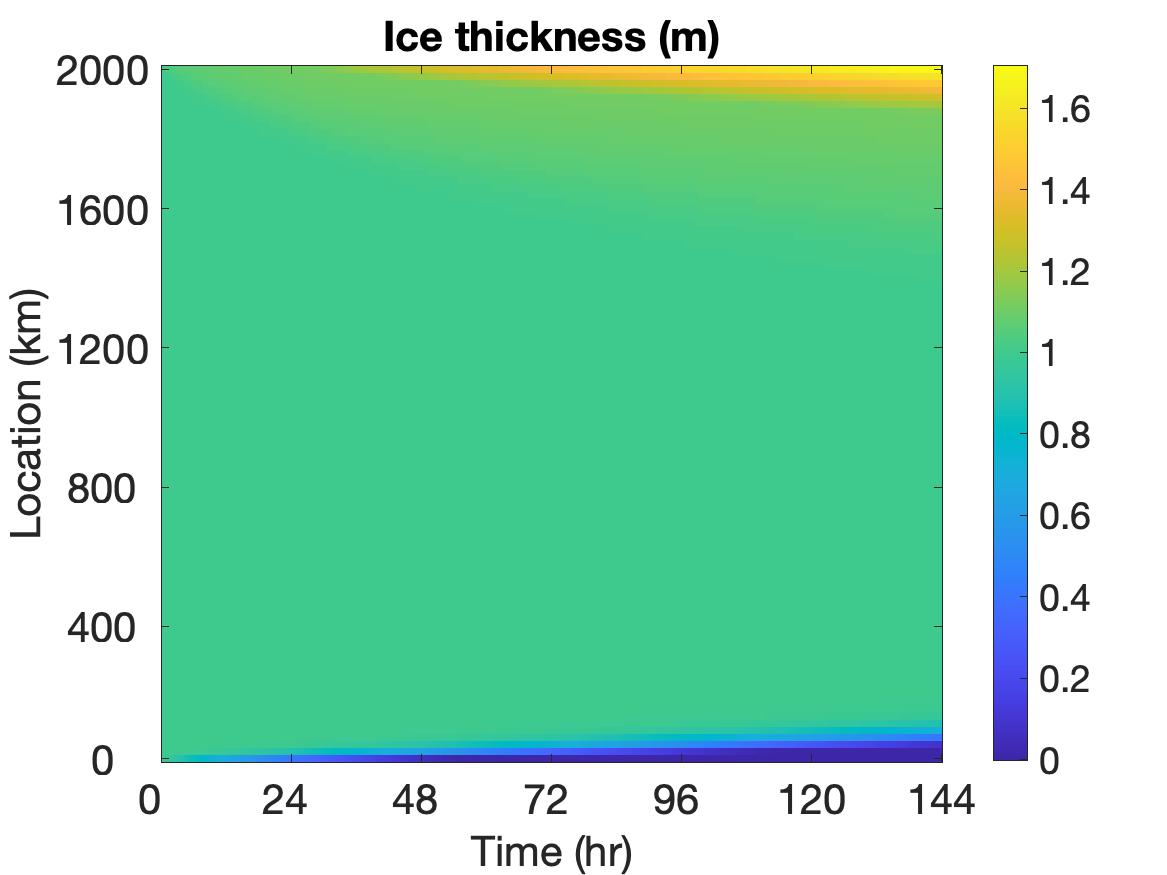}
\includegraphics[width=0.32\textwidth]{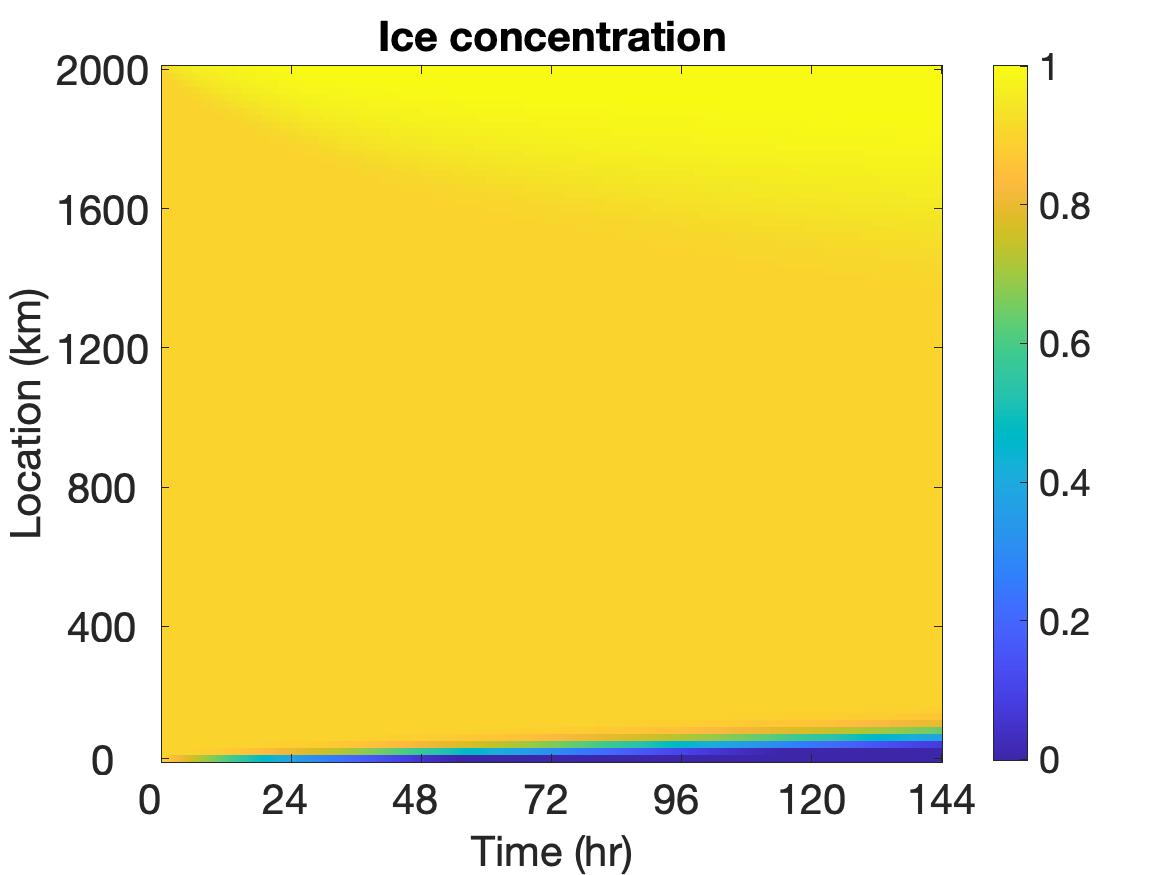}\\[1ex]
\includegraphics[width=0.32\textwidth]{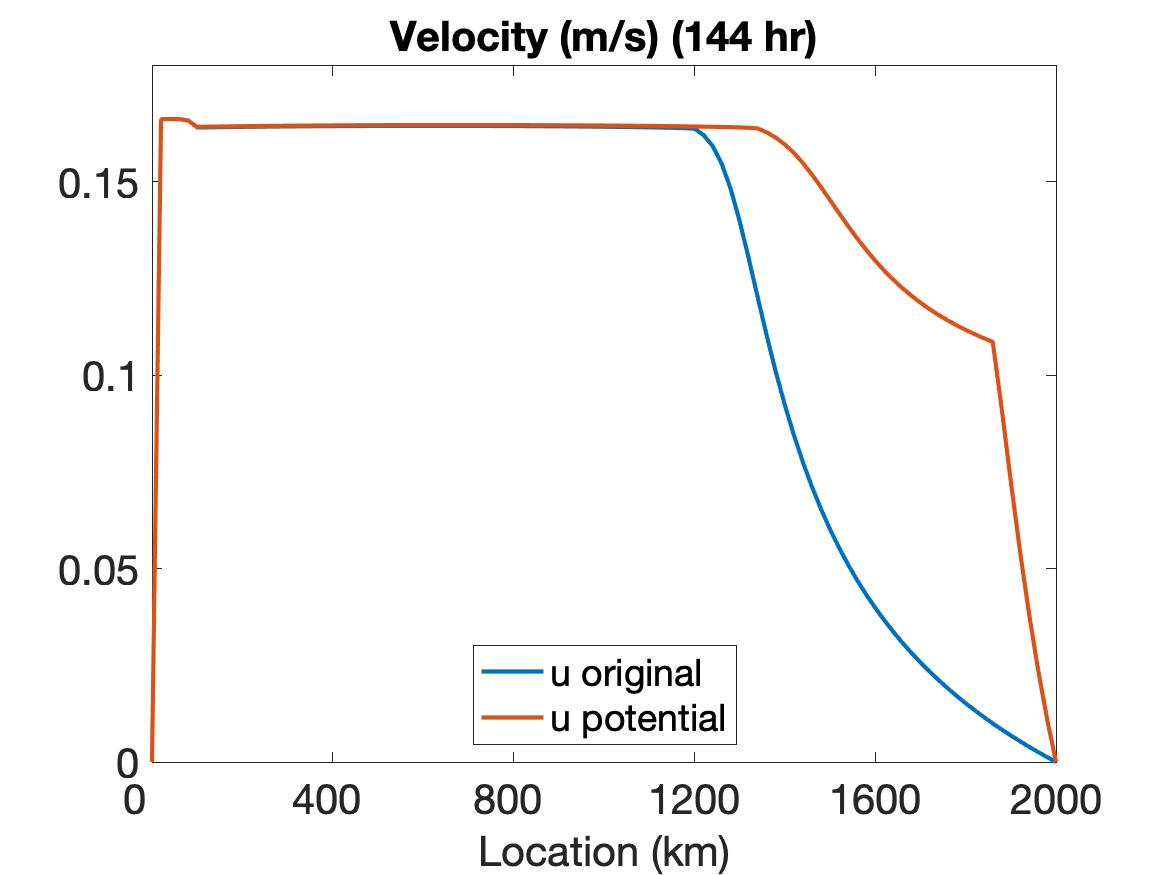}
\includegraphics[width=0.32\textwidth]{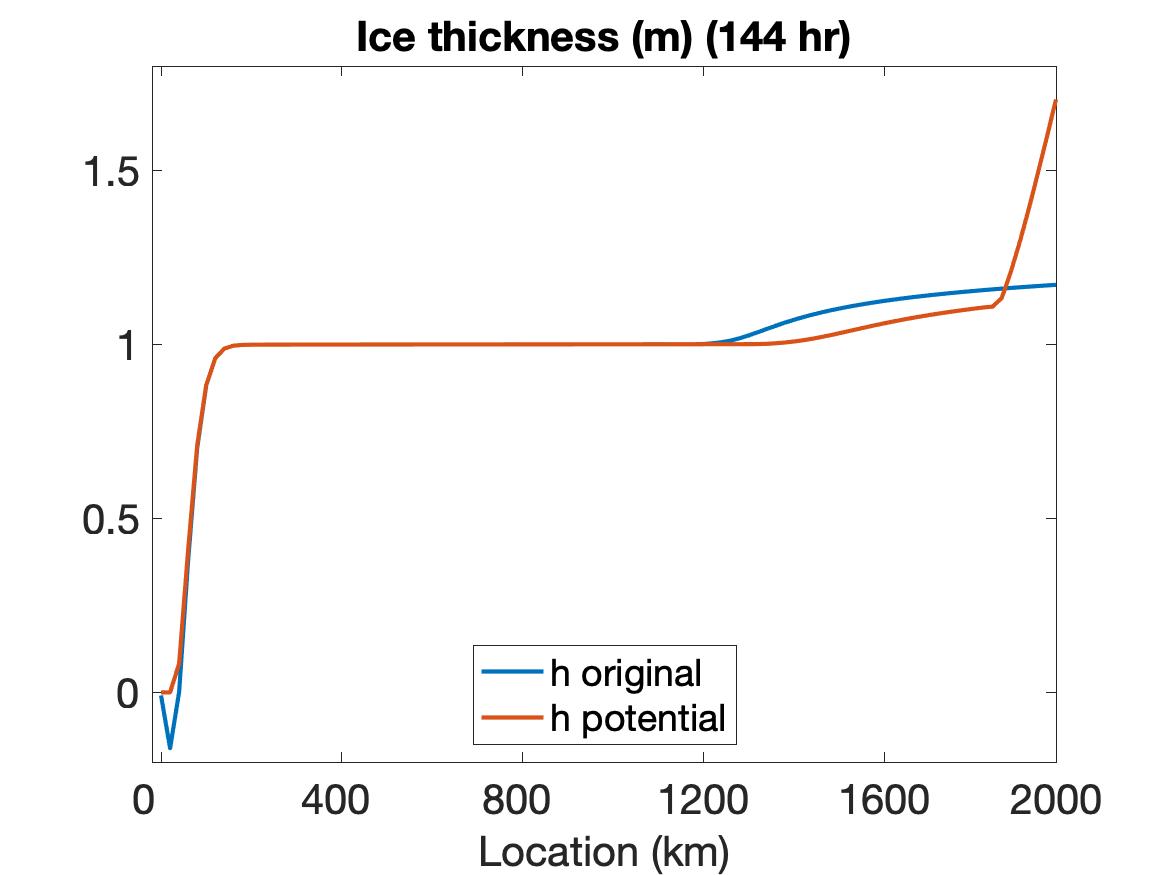}
\includegraphics[width=0.32\textwidth]{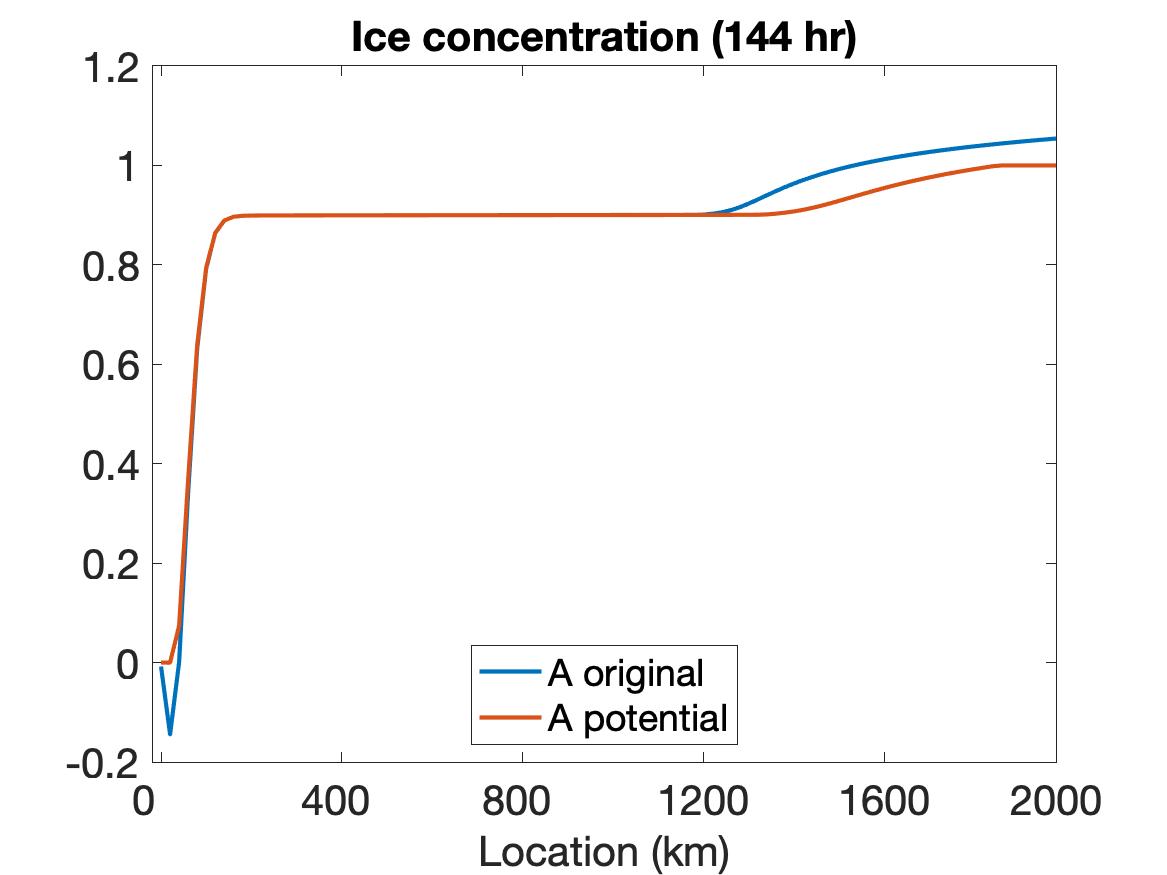}
\end{center}
\caption{Potential function method example. (Top) images of computed solutions without applying potential function method; (middle row) images of computed solutions applying potential function method; (bottom) solutions plots at final time.}
\label{fig: potential}
\end{figure}

The solutions of $u$, $h$, and $A$ are displayed in Figure \ref{fig: potential}, where the top row depicts the solutions without employing the potential function method,  the middle row shows the solutions when the potential function is used, and the bottom row shows the solution plots at the final time.
It is evident that using the potential function yields some differences in the solution. The change of color shade in the velocity image around the 2000 km boundary of the region at roughly 30 hours indicates the situation when the non-physical values are detected and potential function starts to take effect. With the potential function method, the velocity remains a relatively large value on a larger portion of the domain and decreases to 0 in a sharper manner towards the end of the boundary. The behavior of the ice thickness near the boundary, especially around 2000 km, dramatically changes due to the application of the potential function method. In particular, the non-physical negative values near the left boundary are replaced by smoother physically meaningful values. Also, the thickness value is much larger near the right boundary, and the ridging effect is more clearly observed. This makes more sense physically for ridging on the ocean-land boundary. For the ice concentration, non-physical negative values near the left boundary and non-physical large values near the right boundary are also appropriately treated by the potential function method.

\begin{remark}
As shown in the numerical tests, the potential function method plays an important role in preserving the bounds for ice concentration and ice thickness in the sea ice model. Compared with the cut-off approach, which simply removes any values outside the desired range, the potential function method has the advantage of not requiring any post-processing that may introduce discontinuities into the solution. Our results (not shown here) also demonstrate that for the current numerical test case setup, the potential function method and cut-off approach perform similarly with respect to conservation. It remains an open question to see how the potential function method performs in more complicated scenarios, such as in the case of non-uniform wind forcing or in higher-dimensional settings. Such details will be explored in future work.
\end{remark}

\section{Concluding remarks}
\label{sec:conclusion}
This paper discusses the current methodology and limitations, namely poor convergence and out-of-range issues for both ice concentration and ice thickness,  for solving the  VP sea ice model. To improve the performance of the numerical solutions, we propose the use of higher-order methods. In particular, a case study of the celebrated WENO scheme is provided for the one-dimensional sea ice model, and we verify its improved numerical convergence when compared to standardly employed algorithms.  Moreover, WENO is able to resolve discontinuities and sharp features that may occur in sea ice covers.  With regard to the out-of-range issue, this investigation proposes and implements a potential function method that naturally incorporates the physical restrictions of ice thickness and ice concentration in the transport equations.

Since it is relatively easier to examine numerical convergence properties, the current work is restricted to a one-dimensional case study.  Moving forward, we will test the ideas here of using both higher-order methods and the potential function approach in more realistic environments, including the two-dimensional model as well as physical set-up test regimes. Besides viscous-plastic rheology, other rheologies have also been proposed, such as elastic-anisotropic-plastic rheology \cite{Wilchinsky06} and Maxwell elasto-brittle rheology \cite{Dansereau16}. Investigating the numerical performance of the approaches used here may benefit the numerical solutions for these types of rheologies as well. Another avenue for future work is to obtain a more realistic setup of the sharp features in the sea ice cover by incorporating available observations with the physical model via data assimilation techniques. For example, in \cite{Asadi19} data assimilation experiments are based on the one-dimensional VP model discretized by a centered difference scheme, and in future investigations, we can adapt this approach to the higher-order method framework discussed here. Finally, sparse features in the ice thickness have been observed in \cite{Asadi19}, leading to the successful implementation of an $\ell_1-\ell_2$ regularization approach. A general framework for incorporating $\ell_1$ regularization into numerical solvers for partial differential equations with sparse solutions was developed in \cite{Scarnati18}, and combining ideas from there along with the results here may also be beneficial within the data assimilation framework.


\begin{thebibliography}{100}

\bibitem{Arakawa77}{\sc A. Arakawa and V.R. Lamb}, {\it Computational design of the basic dynamical processes of the UCLA general circulation model}. Methods Comput. Phys., 17, 173-265, 1977.

\bibitem{Asadi19}{\sc N. Asadi, K. A. Scott and D. A. Clausi}, {\it Data fusion and data assimilation of ice thickness observations using a regularisation framework}. Tellus A: Dynamic Meteorology and Osceanography, 71:1, 2019.

\bibitem{Auclair17}{\sc J-P. Auclair, J-F. Lemieux, L.B. Tremblay and H. Ritchie}, {\it Implementation of Newton's method with an analytical Jacobian to solve the 1D sea ice momentum equation}. J. Comput. Phys., 340:69-84, 2017.

\bibitem{Dansereau16}{\sc V. Dansereau, J. Weiss, P. Saramito and P. Lattes}, {\it A Maxwell-elasto-brittle rheology for sea ice modelling}. Cryosphere, 10, 1339-1359, 2016.

\bibitem{Fix83}{\sc G.J. Fix}, In {\it Free Boundary Problems: Theory and Applications}, ed. A. Fasano and M. Primicerio, Boston: Piman, 580.

\bibitem{Girard11}{\sc L. Girard, S. Bouillon, J. Weiss, D. Amitrano, T. Fichefet and V. Legat}, {\it A new modeling framework for sea-ice mechanics based on elasto-brittle rheology}. Annals of Glaciology, 52(57), 123-132, 2011.

\bibitem{Hibler79}{\sc W.D. Hibler}, {\it A dynamic thermodynamic sea ice model}. J. Phys. Oceanogr, 9:815-846, 1979.

\bibitem{Hunke01}{\sc E.C. Hunke}, {\it Viscous-plastic sea ice dynamics with the EVP model: Linearization issues}. J. Comput. Phys., 170:18–38, 2001.

\bibitem{Hunke97}{\sc E.C. Hunke and J.K. Dukowicz}, {\it An elastic-viscous-plastic model for sea ice dynamics}. J. Phys. Oceanogr, 27:1849-1867, 1997.

\bibitem{Hunke10}{\sc E.C. Hunke, W. Lipscomb and A. Turner}, {\it Sea-ice models for climate study: Retrospective and new directions}. Journal of Glaciology, 56(200), 1162-1172, 2010.

\bibitem{Hunke15}{\sc E.C. Hunke and W.H. Lipscomb, A.K. Turner, N. Jeffery, S. Elliott}, {\it CICE: the Los Alamos sea ice model documentation and software user's manual version 5.1}. Tech. rep., Los Alamos National Laboratory, 2015.

\bibitem{Ip91}{\sc C.F Ip, W.D. Hibler, G.M. Flato}, {\it On the effect of rheology on seasonal sea-ice simulations}. Ann. Glaciol, 15: 17-25, 1991.

\bibitem{Jiang96}{\sc G. Jiang and C. Shu}, {\it Efficient implementation of weighted ENO schemes}, J. Comput. Phys., 126:202-228, 1996.

\bibitem{Kimmritz15}{\sc M. Kimmritz, S. Danilov, M. Losch}, {\it On the convergence of the modified elastic–viscous–plastic method for solving the sea ice momentum equation}, J. Comput. Phys., 296:90-100, 2015.

\bibitem{Kobayashi10}{\sc R. Kobayashi}, {\it A brief introduction to phase field method}. AIP Conference Proceedings, 1270:282, 2010.

\bibitem{Kuzmin10}{\sc D. Kuzmin}, {\it A Guide to Numerical Methods for Transport Equations}. Friedirch-Alexander Universitäte Erlangen-Nürnberg, 2010.

\bibitem{Kwok08}{\sc R. Kwok, E.C. Hunke, D. Maslowski, and J. Zhang}, {\it Variability of sea ice simulations assessed with RGPS kinematics}. J. Geophys. Res., 113(C11), 2008.

\bibitem{Langer86}{\sc J.S. Langer}, {\it Models of pattern formation in first-order phase transitions}. In {\it Directions in Condensed Matter Physics}, ed. G. Grinstein and G. Mazenko, 165–186. Singapore: World Scientific.

\bibitem{Lemieux08}{\sc J-F. Lemieux, B. Tremblay, S. Thomas, J. Sedl\'a\v cek and L.A. Mysak}, {\it Using the preconditioned Generalized Minimum RESidual(GMRES) method to solve the sea-ice momentum equation}. J. Geophys. Res., 113, 2008.

\bibitem{Lemieux09}{\sc J-F. Lemieux and B. Tremblay}, {\it Numerical convergence of viscous-plastic sea ice models}. J. Geophys. Res., 114, 2009.

\bibitem{Lemieux10}{\sc J.-F. Lemieux, B. Tremblay, J. Sedlácek, P. Tupper, S. Thomas, D. Huard and J.-P. Auclair}, {\it Improving the numerical convergence of viscous-plastic sea ice models with the Jacobian-free Newton Krylov method}. J. Comput. Phys., 229:2840–2852, 2010.

\bibitem{Lemieux12}{\sc J-F. Lemieux, D. Knoll, B. Tremblay, D. Holland and M. Losch}, {\it A comparison of the Jacobian-free Newton-Krylov method and the EVP model for solving the sea ice momentum equation with a viscous-plastic formulation: A serial algorithm study}. J. Comput. Phys., 231:5926–5944, 2012.

\bibitem{Lemieux14}{\sc J.-F. Lemieux, D.A. Knoll, M. Losch and C. Girard}, {\it A second-order accurate in time IMplicit–EXplicit (IMEX) integration scheme for sea ice dynamics}. J. Comput. Phys., 263:375-392, 2014.

\bibitem{Lipscomb04}{\sc W.H. Lipscomb and E.C. Hunke}, {\it Modeling sea ice transport using incremental remapping}. Mon. Weather Rev., 132(6), 1341-1354, 2004.

\bibitem{Lipscomb07}{\sc W.H. Lipscomb, E.C. Hunke, W. Maslowski and J. Jakacki}, {\it Ridging, strength, and stability in high-resolution sea ice models}. J. Geophys. Res., 112, C03S91, 2007.

\bibitem{Liu11}{\sc Y. Liu, C. Shu and M. Zhang}, {\it High order finite difference WENO schemes for nonlinear degenerate parabolic equations}, SIAM J. Sci. Comput., 33(2):939-965, 2011.

\bibitem{Liu94}{\sc X. Liu, S. Osher and T. Chan}, {\it Weighted essentially non-oscillatory schemes}, J. Comput. Phys., 115, 200–212, 1994.

\bibitem{McPhee75}{\sc M. McPhee}, {\it Ice-ocean momentum transfer for the AIDJEX ice model}, A.I.D.J.E.X. Bull., 29:93-111, 1975.

\bibitem{Mehlmann19}{\sc C. Mehlmann}, {\it Efficient numerical methods to solve the viscous-plastic sea ice model at high spatial resolutions}. Dissertation, Otto-von-Guericke Universität Magdeburg, 2019.

\bibitem{Parno19}{\sc M. Parno, B. West, A. Song, T. Hodgdon, DT. O'Conner}, {\it Remote measurement of sea ice dynamics with regularized optimal transport}. Geophysical Research Letters 46(10), 5341–5350, 2019.

\bibitem{Saad93}{\sc Y. Saad}, {\it A flexible inner-outer preconditioned GMRES algorithm}. SIAM J. Sci. Comput, 14:461-469, 1993.

\bibitem{Scarnati18}{\sc T. Scarnati, A. Gelb, R. B. Platte}, {\it Using $\ell_1$ regularization to improve numerical partial differential equation solvers}. J. Sci. Comput., 75(1), 225-252, 2018.

\bibitem{Seinen17}{\sc C. Seinen},
{\it A fast and efficient solver for Viscous-Plastic sea ice dynamics}.
Master’s thesis, University of Victoria, 2017.

\bibitem{Shu20}{\sc C. Shu},
{\it Essentially non-oscillatory and weighted essentially non-oscillatory schemes}. Acta Numerica, 29, 701-762, 2020.

\bibitem{Shu88}{\sc C. Shu and S. Osher},
{\it Efficient implementation of essentially non-oscillatory shock-capturing schemes}. J. Comput. Phys., 77(2), 439-471, 1988.

\bibitem{Wang21}{\sc Q. Wang, L. Ju and Y. Xue}, {\it Chapter 13 - The application of peridynamics for ice modeling}. Editor(s): Erkan Oterkus, Selda Oterkus, Erdogan Madenci,
In Elsevier Series in Mechanics of Advanced Materials, Peridynamic Modeling, Numerical Techniques, and Applications, Elsevier, 275-308, 2021.

\bibitem{Wilchinsky06}{\sc A. V. Wilchinsky and D. L. Feltham}, {\it Modelling the rheology of sea ice as a collection of diamond-shaped floes}. J. Non-Newtonian Fluid Mech., 138, 91-107, 2006.

\bibitem{Williams18}{\sc J. Williams and B. Tremblay}, {\it The dependence of energy dissipation on spatial resolution in a viscous-plastic sea ice model}. Ocean Modelling, 130, 40-47, 2018.

\bibitem{Williams17}{\sc J. Williams, B. Tremblay and J.-F. Lemieux}, {\it The effects of plastic waves on the numerical convergence of the viscous–plastic and elastic–viscous–plastic sea-ice models}. J. Comput. Phys., 340, 519-533, 2017.

\bibitem{Zhang91}{\sc J. Zhang and W.D. Hiber}, {\it On an efficient numerical method for modeling sea ice dynamics}. J. Geophys. Res., 102:8691-8702, 1991.


\end{thebibliography}
\end{document}